\documentclass[a4paper]{amsart}

\usepackage{amssymb,amsmath,graphics,verbatim}
\usepackage{latexsym}
\usepackage{eucal}

\makeatletter
\@addtoreset{equation}{section}\makeatother

\newtheorem{theo}{Theorem}[section]
\newtheorem{lem}[theo]{Lemma}
\newtheorem{prop}[theo]{Proposition}
\newtheorem{cor}[theo]{Corollary}
\theoremstyle{remark} \newtheorem{remark}[theo]{Remark}
\newtheorem{defn}[theo]{Definition}
\newcommand{\mc}{\mathcal}
\newcommand{\rr}{\mathbb{R}}
\newcommand{\nn}{\mathbb{N}}

\newcommand{\zz}{\mathbb{Z}}

\newcommand{\eps}{\epsilon}

\newcommand{\pl}{\partial}
\newcommand{\x}{\times}

\newcommand{\cjd}{\rangle}
\newcommand{\cjg}{\langle}

\newcommand{\demi}{\frac{1}{2}}
\newcommand{\ndemi}{\frac{n}{2}}

\newcommand{\zf}{\textrm{zf}}
\newcommand{\bfo}{\textrm{bf}_0}
\newcommand{\rbo}{\textrm{rb}_0}
\newcommand{\lbo}{\textrm{lb}_0}
\newcommand{\lb}{\textrm{lb}}
\newcommand{\rb}{\textrm{rb}}
\newcommand{\bfa}{\textrm{bf}}
\newcommand{\bfc}{\textrm{bf}}
\newcommand{\sca}{\textrm{sc}}

\newcommand\Id{\operatorname{Id}}

\newcommand\Diff{\textrm{Diff}}
\newcommand\RR{\mathbb{R}}
\newcommand\CC{\mathbb{C}}
\newcommand\NN{\mathbb{N}}
\newcommand\extunion{\overline{\cup}}
\newcommand\Omegab{\tilde \Omega_b^{1/2}}
\newcommand\MMksc{M^2_{k, \sca}}
\newcommand\MMMksc{M^3_{k, \sca}}
\newcommand\phg{\operatorname{phg}}
\newcommand\Omegabht{\tilde \Omega_b^{1/2}}
\renewcommand\Re{\operatorname{Re}}

\newcommand\Mbar{M}

\begin{document}
\title[Resolvent at low energy and Riesz transform]{Resolvent at low energy and Riesz transform 
for Schr\"odinger operators on asymptotically conic manifolds. I.}
\author{Colin Guillarmou}
\address{Laboratoire J. Dieudonn\'e\\
Universit\'e de Nice\\ Parc Valrose\\
06100 Nice\\ FRANCE}     
\email{cguillar@math.unice.fr}
\author{Andrew Hassell}
\address{Department of Mathematics, Australian National University \\ Canberra ACT 0200 \\ AUSTRALIA}
\email{hassell@maths.anu.edu.au}
\subjclass[2000]{}
%

\begin{abstract} Let $M^\circ$ be a complete noncompact manifold and $g$ an asymptotically conic manifold on $M^\circ$, in the sense that  $M^\circ$ compactifies to a manifold with boundary $M$ in such a way that $g$ becomes a scattering metric on $M$. A special case that we focus on is that of asymptotically Euclidean manifolds, where the induced metric at infinity is equal to the standard metric on $S^{n-1}$; such manifolds have an end that can be identified with $\RR^n \setminus B(R,0)$ in such a way that the metric is asymptotic in a precise sense to the flat Euclidean metric. 
We analyze the asymptotics of the resolvent kernel $(P + k^2)^{-1}$ where $P = \Delta_g + V$ is the sum of the positive Laplacian associated to $g$ and a real potential function $V$ that is smooth on $M$ and vanishes to third order at the boundary (i.e. decays to third order at infinity on $M^\circ$). We show that on a blown up version of $M^2 \times [0, k_0]$ the resolvent kernel is conormal to the lifted diagonal and polyhomogeneous at the boundary,
and we are able to identify explicitly the leading behaviour of the kernel at each boundary hypersurface. Using this we show that  the Riesz transform of $P$ is bounded on $L^p(M^\circ)$ for $1 < p < n$, and that this range is optimal if $V \not\equiv 0$ or if $M$ has more than one end. The result with $V\not\equiv 0$ is new even when $M^\circ = \RR^n$, $g$ is the Euclidean metric and $V$ is compactly supported. When $V\equiv 0$ with one end, the range of $p$ becomes $1 < p < p_{\max}$ where 
$p_{\max}>n$ depends explicitly on the first non-zero eigenvalue of the Laplacian on the boundary $\pl M$. 

Our results hold for all dimensions $\geq 3$  under the assumption that $P$ has neither zero modes nor a zero-resonance. In a follow-up paper we shall analyze the same situation in the presence of zero modes and zero-resonances. 
\end{abstract}
\maketitle
\section{Introduction}

We study the resolvent $(P+k^2)^{-1}$ as $k \downarrow 0$ for 
Schr\"odinger operators $P=\Delta_g+V$ on 
asymptotically conic manifolds. Let $M^\circ$ be a smooth manifold with one or more ends, each 
diffeomorphic to $S_i\x(0,\infty)$ where $S_i$ is a smooth compact manifold of dimension $n-1$. Write $S$ for the disjoint union of the $S_i$. 
We suppose $M^\circ$ admits a compactification $\Mbar$ which is smooth with boundary $\pl\Mbar=S$ and
let $g$ be an \emph{asymptotically conic metric} (also called
\emph{scattering metric}) on $M$, that is a complete metric on $M$ such that in a collar neighbourhood 
$[0,\eps)_x\x\pl\Mbar$ near $\pl\Mbar$, 
\begin{equation}\label{metricconic}
g=\frac{dx^2}{x^4}+\frac{h(x)}{x^2}
\end{equation}
where $x$ is a smooth function that defines the boundary $\pl\Mbar$ (i.e. $\pl\Mbar=\{x=0\}$ and $dx$ does not vanish on $\pl\Mbar$) and $h(x)$ is a smooth family of metrics on $S$. 
We let $V$ be a real potential such that 
\begin{equation}\label{hyp2}
V\in C^{\infty}(\Mbar), \quad V=O(x^{3}) \textrm{ as }x\to 0 .
\end{equation}
A particular case which will be of special interest in this work is the case where $S=S^{n-1}$ is the sphere and
$h_0:=h(0)=d\theta^2$ is its canonical metric. Then we say that $g$ is \emph{asymptotically Euclidean} since the
Euclidean metric in polar coordinates reads $dr^2+r^2d\theta^2$ on $\rr_{r}^+\x S^{n-1}$ and can be put under the form
$dx^2/x^{4}+d\theta^2/x^{2}$ by setting $x:=1/r$. For a general scattering metric taking the form \eqref{metricconic}, we can view $r = 1/x$ as a generalized `radial coordinate', as the distance to any fixed point of $M$ is given by $r + O(1)$ as $r \to \infty$. 

Let $\Delta_g$ be the positive Laplacian on $M$ with respect to $g$. 
Then the Schr\"odinger operator   
\[P=\Delta_g+V\]
is self-adjoint on $L^2(M,dg)$ and its spectrum is $\sigma(P)=[0,\infty)\cup \sigma_{\rm pp}(P)$
where $\sigma_{\rm pp}(P)=\{-k^2_1\geq \dots \geq-k^2_N\}$ is a finite set of negative eigenvalues (by convention $k_i>0$). 
Note that $0$ can
be an $L^2$-eigenvalue but the half-line $(0,\infty)$ is only continuous spectrum.\\ 

\begin{figure}[ht!]\label{mmksc}
\begin{center}
\input{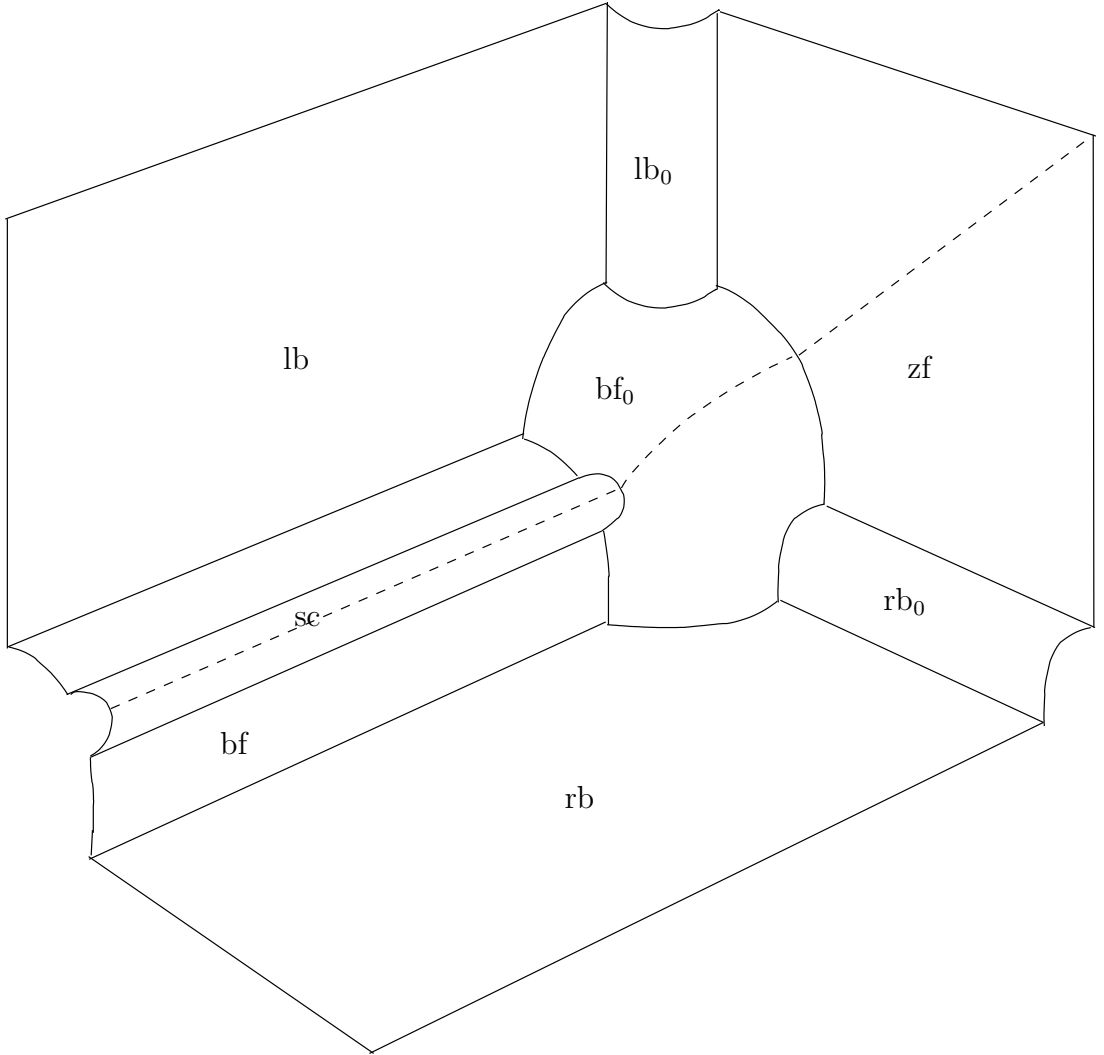}
\caption{The blow-up manifold $M^2_{k,\sca}$; the dashed line is the boundary of $\Delta_{k, \sca}$}
\end{center}
\end{figure}

The resolvent $R(k)=(P+k^2)^{-1}$ is well defined as a bounded operator on $L^2(M,dg)$
for $k \in (0, k_1)$ but fails to be bounded or defined at $k=0$. We describe in this work 
the uniform behaviour near $k=0$ in terms of distributional kernel of $R(k;m,m')$
on the manifold with corners $X=[0,k_0]_k\x \Mbar\x \Mbar$ for small $k_0 < k_1$. More precisely, we follow Melrose's program \cite{Kyoto}, pass to a desingularized version of $X$,  
that we call $M^2_{k,\sca}$ (see figure~\ref{mmksc}), obtained by blowing up a certain number of corners of $X$, and lift the kernel of $R(k)$ to $\MMksc$. Our result is that $R(k)$ is conormal to the lifted diagonal of $\MMksc$ and polyhomogeneous conormal at the boundary. Thus, 
when lifted to $M^2_{k,\sca}$ the kernel has asymptotic expansions at all orders at each face. 
In particular, the expansion at the face corresponding to $M^2 \times \{ k=0 \}$ is simply the Laurent
expansion of $R(k)$ at $k=0$. 

This Laurent expansion is quite important if one wants to
understand the large time asymptotics of Schr\"odinger propagators, heat operators or
wave operators. It has been studied by Jensen-Kato \cite{JK}
for $(M^\circ,g)=(\rr^3, g_{\rm eucl})$ and extended by Murata \cite{Mu} for $M^\circ=\rr^n$ and $P=p(D)+V$
a real constant coefficient elliptic differential operator with $V$ a compact operator in 
weighted Sobolev spaces. Then more recently X.P. Wang \cite{XPW} studied the case of the Laplacian plus a potential which is homogeneous of degree $-2$ on an exact conic 
manifold,  perturbed by a shorter range potential. He showed there is an asymptotic expansion of $R(k)$
at $k=0$ with a number of terms that depend on the decay rate of the perturbation.  
In all these approaches, an exact formula for the free model operator
is used to derive the Laurent expansion and compute its leading coefficients, which are expressed
in terms of the model resolvent. It is not clear how to generalize this approach to cases with no free model. 
  
Our approach is to constuct the resolvent kernel through a parametrix, guessing what must be the 
Taylor coefficients at each face of $M^2_{k,\sca}$ and in particular at $k=0$. This, in turn, relates 
to the analysis of a Laplacian type operator on a manifold $(M,g_b)$ conformal to $(M,g)$,
\[g_b:=x^2g=\frac{dx^2}{x^2}+h(x),\]
that is, a manifold with asymptotically \emph{cylindrical} ends. We use the b-calculus of \cite{APS}
to solve the model problems at $k=0$  and, as a corollary, we  describe the first coefficients
of the Laurent expansion of $R(k)$ at $k=0$. Our results are phrased in terms of a calculus of pseudodifferential operators $\Psi^{m, (a_{\bfo}, a_{\zf}, a_{\sca}); \mc{E}}_k(M; \Omegab)$ acting on half-densities on $M$, where $m$ is the pseudodifferential order, $(a_{\bfo}, a_{\zf}, a_{\sca})$ determines the growth of the Schwartz kernel of the operator at the three boundary hypersurfaces meeting the diagonal (see Figure~\ref{mmksc}) and $\mc{E}$ is an `index family' listing the allowable terms in the asymptotic expansion of the kernel at each boundary hypersurface of $\MMksc$. 

\begin{theo} Suppose that $n \geq 3$ and that  $P$ has neither zero-modes nor a zero-resonance\footnote{By a zero mode we mean an $L^2$ eigenfunction with eigenvalue $0$, and by a zero-resonance a solution $v$ to $Pv = 0$ that decays to zero at infinity. Zero-resonances can only exist for $n \leq 4$.}. Then 
\begin{equation}
R(k) \in \Psi^{-2, (-2, 0, 0), \mc{R}}(M, \Omegab)
\label{res-kernel}\end{equation}
where $\mc{R}$ is an index set. This index set is given explicitly in \eqref{rif} when $M$ is asymptotically Euclidean and lower bounds are given in \eqref{rif2} for the general case. More informally, $R(k)$ is conormal to the diagonal of $M^2_{k,\sca}$ and polyhomogeneous at the boundary, uniformly down to $k=0$. 
\end{theo}

We use this result to analyze the boundedness of the so-called 
Riesz transform as an operator from $L^p(M,dg)$ to $L^p(M,T^*M,dg)$. If $\Pi_>$ denotes the spectral projection for the operator $P$ onto the positive spectrum $(0, \infty)$, and $P_> = P \circ \Pi_>$ then we define the Riesz transform $T$ by 
\[ T:= d \circ (P_>)^{-\demi}.
\]
Here $d$ is exterior differentiation from functions to $1$-forms.\\

The Riesz transform for the Laplacian on a complete non-compact manifold is clearly bounded on $L^2$ 
but its $L^p$ boundedness is in general quite involved and can be interpreted as a 
way of comparing two definitions of Sobolev spaces. We refer the reader to the paper \cite{ACDH}
for general results and its section 1.3 for the state of the art about Riesz transform $L^p$ boundedness.
Here we mention that for the Laplacian, $T\in\mc{L}(L^p)$ if $p\in(1,2]$ 
in a quite general setting (see Coulhon-Duong \cite{CD}); in other words low ranges of $p$
seem less sensitive to the geometry. For $p>2$, the question is more delicate
and the first cases to look at are manifolds which have a Euclidean structure near infinity.
Using the formula
\[
P_>^{-\demi}=\frac{2}{\pi}  \int_0^\infty R(k) \Pi_> \, dk
\]
Carron-Coulhon-Hassell \cite{CCH} analyzed the $L^p$ boundedness properties of $T$ in this setting 
(again $V=0$) through a careful analysis of the resolvent. They showed
that $T\in\mc{L}(L^p)$ if and only if $p\in(1,\infty)$ when the manifold has one end, whereas $p\in(1,n)$ 
when there is more than one end. The unboundedness for $p > n$
had been shown earlier by Coulhon-Duong \cite{CD}. 
 
As for the case $V\not\equiv 0$, some work has been done on $L^p$ estimates, 
especially for Schr\"odinger operators in $\rr^n$,
but in most cases, it seems that the potential is supposed \emph{non-negative} so heat kernel techniques are effective. 
We note for instance the boundedness of Riesz transform on $L^p$ for $p\in(1,2]$ if $V\geq 0$
and $V\in L^1_{\rm loc}(\rr^n)$ (see E.M. Ouhabaz \cite{EMO}), which means again a quite strong stability
for low range of $p$. A similar result is obtained by Coulhon-Duong \cite{CD} for the Laplacian 
on a class of manifolds. P. Auscher and A. Ben Ali \cite{ABA} 
improved a result of Z. Shen \cite{SH} on higher ranges of $p$, this is when the potential 
is non-negative and satisfy the $q$ reverse H\"older inequality for some $q\in(1,+\infty]$, 
indeed they prove in this setting that there is $\eps>0$ such that 
$T\in \mc{L}(L^p)$ if $p\in(1,q+\eps)$ (Z. Shen's restriction was $n/2\leq q<n$).  
We remark that powers $|z|^{-\alpha}$ satisfy the reverse H\"older inequality if $\alpha\in(-\infty,n/q)$. It is interesting that their hypothesis rules out potentials $V$ decaying too fast near infinity, while our hypothesis by contrast requires sufficiently fast decay at infinity. 

Continuing the program of Carron-Coulhon-Hassell \cite{CCH} (see their last section), 
we extend their results to more general manifolds and include the case with 
potential, which need not be non-negative.
We do assume in this paper that $P$ has neither zero modes, nor a zero-resonance. We shall see, in a sequel to this paper, that the presence of zero modes and zero-resonances has dramatic consequences for the Riesz transform 
$L^p$ boundedness. This analysis in turn should facilitate the analysis of the Riesz transform for differential forms on asymptotically conic manifolds, since then $L^2$ harmonic forms and resonances may appear.

We first state our Riesz transform results in the case $V=0$.
\begin{theo}\label{v=0}
Let $n \geq 3$, and let $P$ be the Laplacian on an asymptotically conic manifold $(M,g)$ with one end.
Let $\Delta_{\pl\Mbar}$ be the Laplacian on the boundary of $\Mbar$ for the 
metric $h(0)$ given by (\ref{metricconic}),  let $\lambda_1$ be its first non-zero eigenvalue, and let $\nu_1 = \sqrt{((n-2)/2)^2 + \lambda_1}$. 
If $\nu_1<n/2$, then the Riesz transform $T$ satisfies 
\[ T\in \mc{L}(L^p(M),L^p(M;T^*M)) \iff 1<p<n \Big(\frac{n}{2} - \nu_1\Big)^{-1} \] 
while if $\nu_1\geq n/2$, then
 \[ T\in \mc{L}(L^p(M),L^p(M;T^*M)),\quad \textrm{ for all }\quad  1<p<\infty.\] 
\end{theo}
Notice that $n (n/2 - \nu_1)^{-1}$, the upper threshold for $p$ when $\nu_1 < n/2$,  is always greater than $n$. 
This result answers Problem 8.1 of \cite{CCH}. It extends that of \cite{CCH} for Euclidean ends and is closely related to (and indeed uses) that of Li \cite{Li}, who proves $L^p$ boundedness in the same range for the Riesz transform on metric cones. One direction of these equivalences can also be obtained by the recent paper of Coulhon-Dungey \cite{CD} with  a different method. 
Next, we show 
\begin{theo}\label{th2}
Let $n\geq 3$ and let $P=\Delta_g+V$ with $V$ satisfying  (\ref{hyp2})  
 and such that  $P$ has no zero modes or zero-resonance. 
Suppose that either $M$ has more than one end, or that 
$V\not\equiv 0$. Then
\[ T\in \mc{L}(L^p(M),L^p(M;T^*M)) \iff 1<p<n.\]
\end{theo}

\begin{remark} This seems surprising at first, since $P$ may be a very mild perturbation of the flat Laplacian on $\RR^n$ (if $V$ is, say, compactly supported and nonnegative). Yet it reduces  the upper threshold for $L^p$ boundedness of $T$  from $\infty$ to $n$. The lack of boundedness for the Riesz transform for $p > n$ for the case of the connected sum of two copies of $\RR^n$ was remarked on by Coulhon and Duong \cite{CD}, but there it appeared to be the global geometry and topology of $M$ that was responsible for this phenomenon. However, \cite{CCH} and the present paper give a different perspective; it is \emph{whether there is a bounded nonconstant solution to $Pv = 0$} that determines whether the Riesz transform is unbounded for $p > n$. 
\end{remark}

\begin{remark}
To obtain these Riesz transform results it is essential that we understand the behaviour of $R(k)$ at \emph{all} the boundary hypersurfaces at $k=0$, and not just at the lift of $M^2 \times \{ 0 \}$. Neither Jensen-Kato nor Murata do this, and consequently, it does not seem possible to deduce boundedness properties of the Riesz transform from their results. 
\end{remark} 

\begin{remark}
It can be seen from the Hardy inequality that if  
$M=\rr^n$ with $n>2$ and such that $|V_-(z)|\leq \alpha|z|^{-2}$ 
for $\alpha<(n/2-1)^2$ where $V_-$ is the negative part of $V$, then $P=\Delta_{\rr^n}+V$
does not have  a zero-resonance nor nontrivial $L^2$ kernel. In particular, the 
frequently-made assumption $V_-=0$  avoids such problems.
\end{remark}

\emph{Acknowledgements.} Much of this work was done during a visit of  C.G. to the Mathematics Department, Australian National University funded by 
Australian Research Council Discovery Grant DP0449901.  
The research was supported in part also by NSF 
grant DMS\-0500788 and French ANR grants JC05-52556 and JC0546063 (C. G.). 
The authors thank Richard Melrose and Ant\^onio S\'a Barreto for kindly supplying them with a copy of the unpublished note \cite{Mel}. They also thank  Gilles Carron, Thierry Coulhon, Alan Mc Intosh, Rafe Mazzeo and Adam Sikora for helpful conversations. 

\section{Preliminaries}

\subsection{Manifolds with cylindrical ends}
As mentioned in the introduction, we convert the problem of inverting $P + k^2$ to an equivalent problem on a Riemannian manifold with asymptotically cylindrical ends. Let us assume that $M$ has only one end $E\simeq \rr_r^+\x\pl\Mbar$ 
to simplify the exposition.
When the manifold is asymptotically Euclidean, we have coordinates $z=(ry_1,\dots,ry_n)$ on $E$ induced by the coordinates
$y\in\rr^n$ on the sphere $\{|y|=1\}$, and we denote by $x=|z|^{-1}=1/r$ a smooth function that defines the boundary of the radial compactification $\Mbar$.
We denote by $\mc{V}_b$ the Lie algebra of smooth vector fields on $\Mbar$ that are tangent to the 
boundary $\pl\Mbar$ and $\textrm{Diff}_b(\Mbar)$ the enveloping algebra of $\mc{V}_b$ over $C^\infty(\Mbar)$, 
i.e. the set of differential operators in these vector fields. Such vector fields can be locally written
as a combination ovr $C^\infty(\Mbar)$ of $x\pl_x,\pl_{y_1},\dots,\pl_{y_{n-1}}$ where $(\pl_{y_i})_i$
form a local basis of the tangent space of the boundary $\pl\Mbar$. We denote by $\mc{V}_{\rm sc}$  
the set of smooth vector fields which can locally be written near the boundary as a combination of
$x^2\pl_x, x\pl_{y_1},\dots,x\pl_{y_{n-1}}$ and 
$\textrm{Diff}_{\rm sc}(\Mbar)$ the enveloping algebra of $\mc{V}_{\rm sc}$.\\

We define $g_b:=x^2g$ the conformal metric to $g$ and $L^2_b(M)=L^2(M,{\rm dvol}_{g_b})=x^{-\ndemi}L^2(M)$. 
$(M,g_b)$ is then an exact b-metric in the sense of Melrose \cite{APS}, that is a cylindrical metric on $M$. We now define the operator $P_b$ by 
\begin{equation}\label{pvspb}
P=x^{\ndemi+1}P_bx^{-\ndemi+1}.
\end{equation} 
Since $P$ is formally self-adjoint with respect to $g$, $P_b$ is formally self-adjoint with respect to $g_b$. A calculation gives
\begin{equation}\label{pb-form}
P_b:=-(x\pl_x)^2+\Big(\frac{n-2}{2}\Big)^2+\Delta_{\pl\Mbar}+
W,
\end{equation}
with $W\in x\textrm{Diff}_b(\bar{X})$ a lower-order term at $x=0$.

We study $P$ indirectly via analyzing the operator $P_b$ using the b-calculus and results from \cite{APS}. We now summarize the major results of the b-calculus  as they apply to our particular operator $P_b$. First we recall the b-Sobolev spaces $H^j_b(M)$, $j = 0, 1, 2, \dots$, consisting of those functions in $L^2_b(M)$ which are mapped into $L^2_b(M)$ by all order $j$ elements of $\Diff_b(M)$. The weighted b-Sobolev space $x^\alpha H^j_b(M)$ is defined as\footnote{The b-Sobolev spaces can be defined for any real $j$, but it is not necessary to do so for our purposes.} the space of functions which can be written $x^\alpha f$ with $f \in H^j_b(M)$. 
Let $\lambda_0 = 0 < \lambda_1 \leq \lambda_2 \dots$ be the spectrum of the Laplacian on $(\pl\Mbar,h_0)$. We then define 
\begin{equation}\label{npl}
N_{\pl}:=\{ \nu_0, \nu_1, \dots \mid \nu_i = \sqrt{((n-2)/2)^2 + \lambda_i  } \}.\end{equation}
In the case of the canonical sphere $(\pl\Mbar=S^{n-1},h_0=d\theta^2)$ we have 
\[N_{\pl}=\{(n-2)/2+l; l\in\nn_0\}.
 \]
Melrose's Relative Index Theorem (\cite{APS}, section 6.2) implies in our setting
\begin{theo}\label{rit} The operator $P_b$ is Fredholm as a map from $x^\alpha H^j_b(M)$  to $x^\alpha H^{j-2}_b(M)$ for all $j \geq 2$ and all $\alpha \neq \pm \nu_l$, $l = 0, 1, 2, \dots$. The index of $P_b$ is equal to $0$ for $|\alpha| < (n-2)/2$ and the index jumps by $d_{l}$, the multiplicity of the $\nu_l^2$-eigenspace $E_{\nu_l}$ of $\Delta_{\pl\Mbar}+(n-2)^2/4$, as $\alpha$ crosses the value $\pm\nu_l$ (with $\alpha$ decreasing).  
\end{theo}

We also need results on the regularity of solutions of the equation $P_b u = f$, when $f$ is polyhomogeneous. The following result proved in \cite{APS} is phrased in terms of index sets and the operation of extended union; see subsection~\ref{pfis}. 

\begin{theo}\label{reg} Suppose that $f$ is polyhomogeneous on $M$ with respect to the index set $E$, that $u \in x^\alpha L^2_b(M)$, and that $P_b u = f$. Let $S$ be the set
$$S = \{ \pm \nu_l \mid l = 0, 1, 2, \dots \}
$$
and for $b \in \RR$ let $\mu(b,j)=\sharp \{b+k;k=0,\dots j\}\cap S$. 
Then $u$ is polyhomogeneous with respect to the index set $E \extunion F$, where $F$ is the index set 
$$
\big\{ \big((\pm\nu_l+j), k\big) \mid  l\in\nn_0, j\in\nn_0, \ \pm \nu_l\geq \alpha, \ 0 \leq k \leq \mu(\pm\nu_l,j) - 1 \big\}.
$$
When $(\pl\Mbar,h_0)=(S^{n-1},d\theta^2)$, this reduces to 
\[\big\{ \big(n/2+l, k\big) \mid  l\in\zz, \ n/2+l\geq \alpha, \ 0 \leq k \leq N_l - 1 \big\}.\]
where $N_l$ is the number of elements of the form $\pm(n/2-1+j), j\in\nn_0$ in the interval $(\alpha,n/2+l]$ .
\end{theo}

As a consequence of these two theorems, the vector space
$$
\{ f \in x^{\nu_l - \epsilon} H^j_b(M) \mid P_b f = 0 \} \big/ \{  f \in x^{\nu_l + \epsilon} H^j_b(M) \mid P_b f = 0 \} 
$$
for sufficiently small $\epsilon > 0$ is finite dimensional, independent of $\epsilon$, and  independent of $j$ by elliptic regularity. By Theorem~\ref{reg}, elements of this space have the form $x^{\nu_l} \phi_l(y) + O(x^{\nu_l + \epsilon})$ where $\phi_l(y)$ is an element of $E_{\nu_l}$.  The span of such elements $\phi_l$ is a vector subspace $G_{\nu_l} \subset E_{|\nu_l|}$.

\begin{prop}\label{prop:complementary} (\cite{APS}, Chapter 6) The subspaces $G_{\nu_l}$ and $G_{-\nu_l}$ of $E_{\nu_l}$ are orthogonal complements with respect to the inner product on $(\pl\Mbar,h_0)$. 
\end{prop}

\begin{remark} Thus, the index jumps by $d_{l}$ as $\alpha$ crosses the value $\nu_l$ since the dimension of the null space increases by dim $G_{\nu_l}$, while the dimension of the cokernel,  which can be regarded as the dimension of the null space of $P_b$ on $x^{-l} H^k_b(M)$, increases by dim $G_{-\nu_l}$. 
\end{remark}

\subsection{The blown-up space $M^2_{k,\sca}$}
\subsubsection{Definition of $\MMksc$}

We want to analyze the distributional Schwartz kernel of 
\[R(k):=(P+k^2)^{-1}\]  
as $k\to 0$ and $k>0$.
The model case with $X=\rr^3, g=g_{\rm eucl}$ and $V=0$ gives 
\[R(k;z,z')= \frac1{4\pi} \frac{e^{-k|z-z'|}}{|z-z'|}\] 
this shows that the different asymptotic behaviours of this kernel on $[0,1]_k\x \Mbar\x \Mbar$ leads to consider it
on a blow-up space $M^2_{k,\sca}$ of $[0,1]\x \Mbar\x\Mbar$ as described in a note of Melrose-S\'a Barreto \cite{Mel}. 
For instance $R(k;z,z')$ goes from decreasing at all orders (i.e. $O(x^{\infty})$) as $x=|z|^{-1}\to 0$ when $k>0,z\not=z'$ to decreasing at order $1$ (i.e. $O(x)$) when $k=0$. On the other hand,  in the asymptotic regime $(x,k)\to 0$, $z'$ fixed in $M$, the kernel is a smooth function of the variables $z', z/|z|$, $x$ and $k/x$. This suggests that `polar coordinates' in $(x,k)$ are more adapted to the analysis of $R(k,z,z')$. We realize these and other singular coordinates geometrically by working on the Melrose-S\'a Barreto space in which they appear as smooth functions. 

We use the unpublished note \cite{Mel} and  refer to \cite{cocdmc}, \cite{scatmet} and \cite{APS} for information and notation about manifold with corners, conormal distributions, blow-ups and so on. Here we recall that a submanifold $Y$ of a manifold with corners $X$ 
is said to be a product-type submanifold  or p-submanifold if for every $y \in Y$ there exist local coordinates $x_1, \dots, x_r, y', y''$ with $x_i \geq 0$ local boundary defining functions and $y' \in \RR^{s'}$, $y'' \in \RR^{s''}$, $r + s' + s'' = \dim X$ so that $Y$ is locally defined by  $\{ x_1 = \dots =  x_q = 0, y' = 0 \}$ ($q \leq r$). 
We denote by $[X;Y_1,\dots,Y_N]$ the iterated real blow-up of $X$ around 
$N$ submanifolds $Y_i$ if $Y_1$ is a p-submanifold, the lift of $Y_2$ to $[X; Y_1]$ is a p-submanifold, and so on. 
We shall denote by $\rho_H$ a general boundary defining function for a boundary hypersurface $H$ of $X$.

We now define the space $\MMksc$. Consider in $[0,\eps)\x\Mbar\x\Mbar$ the codimension 3 corner $C_3:=\{0\}\x\pl\Mbar\x\pl\Mbar$ 
and the codimension $2$ corners 
\[C_{2,L}:=\{0\}\x\pl\Mbar\x \Mbar, \quad C_{2,R}:=\{0\}\x\Mbar\x\pl\Mbar,\quad C_{2,C}:=[0,1]\x\pl\Mbar\x\pl\Mbar.\]
We consider first the blow-up 
\[M_{k,b}^2:=[[0,1]\x\Mbar\x\Mbar; C_3, C_{2,R},C_{2,L},C_{2,C}]\]
with blow-down map $\beta_b:M^2_{k,b}\to [0,1]\x\Mbar\x\Mbar$. 
We have $7$ faces on $M^2_{k,b}$, the right, left, and zero faces
\[\rb=\textrm{clos }\beta_b^{-1}([0,1]\x M\x\pl\Mbar),\quad 
\lb:=\textrm{clos }\beta_b^{-1}([0,1]\x\pl\Mbar\x M),\] 
\[\zf:=\textrm{clos }\beta_b^{-1}(\{0\}\x M\x M),\]
the big face $\bfo:=\beta^{-1}(C_3)$ and the b-faces
\[\rbo:=\textrm{clos }\beta_b^{-1}(C_{2,R}\setminus C_3), \quad 
\lbo:=\textrm{clos }\beta_b^{-1}(C_{2,L}\setminus C_3), 
\quad \bfc:=\textrm{clos }\beta^{-1}(C_{2,C}\setminus C_3).\]
The closed lifted diagonal $\Delta_{k,b}=\textrm{clos }\beta_b^{-1}([0,1]\x \{(m,m);m\in M\})$
intersects the face $\bfc$ in a p-submanifold denoted $\pl_{\bfc}\Delta_{k,b}$. We then define
the final blow-up
\[M^2_{k,\sca}:=[M^2_{k,b}; \pl_{\bfc}\Delta_{k,b}].\]

We obtain a set of faces
denoted $\mc{F}=\{\zf,\bfo,\rbo,\lbo,\rb,\lb,\bfa,\sca\}$ by lifting the boundary hypersurfaces of $M^2_b$ 
and denoting $\sca$ the new face obtained by the last blow-up.
These are shown on Figure 1 and we denote by 
$\beta:M^2_{k,\sca}\to M\x M \x[0,1]_k$ the blow-down map, and $\rho_{\rm f}$ a boundary defining function of 
${\rm f}$ for ${\rm f}\in \mc{F}$.
The lifted diagonal on $M^2_{k,\sca}$ is defined by
\[\Delta_{k, \sca}:=\textrm{clos }\beta^{-1}((0,1]\x\{(m,m),m\in M\}).\]

\begin{remark}\label{essential} The essential features of the space $\MMksc$ are
\begin{itemize}
\item for each $k_1 \in (0, k_0)$ the manifold $\MMksc \cap \{ k = k_1\}$ is canonically diffeomorphic to $M^2_{\sca}$;

\item The face $\zf$ is canonically diffeomorphic to $M^2_b$;

\item The  lifted diagonal $\Delta_{k, \sca}$ is a p-submanifold; 

\item The vector fields $\partial_{z_i}$ lifted from the left or the right factor of $M$ to $\MMksc$ are tangent to the boundary and transverse to $\Delta_{k, \sca}$ except at $\bfo$ where they vanish simply; and

\item The face $\zf$ is disjoint from $\sca$, reflecting the different types of asymptotics exhibited by the resolvent for $k =0$ and $k> 0$. The intermediate face $\bfo$ is a transition region between the two types of asymptotics. 
\end{itemize}
\end{remark}

\begin{remark} This space is closely related to, but not the same as, the space in \cite{CCH}. Either space works for asymptotically Euclidean manifolds, but the Melrose-S\'a Barreto space used here seems to be the `correct' one
for asymptotically conic manifolds. 
\end{remark}

\subsubsection{Densities and half-densities}
On any smooth manifold with corners $X$, we denote $\Omega^\demi(X)$ the bundle of smooth half-densities,
$\Omega_{b}^\demi(X)$ the bundle of b-half densities. That is,  $\Omega_{b}^\demi(X)=\rho^{-\demi}\Omega^\demi(X)$ 
if $\rho$ is a product of boundary defining functions for each boundary hypersurface  of $X$. For each boundary hypersurface $H$ of $X$ there is a canonical isomorphism between  $\Omega_{b}^\demi(X) \, | \, H$ and $\Omega_b(H)$, given by cancelling a $|d\rho_H/\rho_H|^\demi$ factor at $\rho_H = 0$, where $\rho_H$ is any  boundary defining function for $H$. This does not depend on the choice of $\rho_H$.   Correspondingly, a smooth b-half-density restricts canonically to a b-half-density on each boundary face $H$. 

We shall denote by $\Omegabht(\MMksc)$ the lift of $\Omegab(M^2_{k,b})$ to $\MMksc$. A smooth nowhere vanishing section of this bundle is $\nu = |dg_b dg'_b dk/k|^\demi$. 
Notice that this bundle is not the b-half-density bundle of $\MMksc$. 
It can be identified with the b-half-density bundle except near the face $\sca$. Indeed, as a b-half-density $\nu$ vanishes to order $n/2$ at $\sca$. We shall express our operators in terms of the bundle $\Omegabht(\MMksc)$ rather than $\Omegab(\MMksc)$, precisely because is it generated by $\nu$ (in the sense that $\nu$ is a smooth nonvanishing section of this bundle).

\subsubsection{Local coordinates and differential operators}
Let us first consider the kernel of the identity operator acting on half-densities. After multiplying by the formal factor $|dk/k|^{1/2}$, this may be written 
\begin{equation}\label{kernelid}
\delta(x-x')\delta(y-y')\left|\frac{dxdx'dydy'dk}{k}\right|^{\demi}.
\end{equation}
Let us write this in coordinates adapted to $\MMksc$ near $\zf$, $\bfo$ and $\sca$. The coordinates above are valid near $\zf$ provided we stay away from $\bfo$.  Near $\zf \cap \bfo$ and in a neighbourhood of $\Delta_{k, \sca}$ we can use $x, \sigma = x'/x, y, y'$ and $\kappa' = k/x'$. In terms of these we can write the identity as
$$
\delta(\sigma - 1) \delta(y-y') \left|\frac{dxd\sigma dydy'd\kappa'}{x  \kappa'}\right|^{\demi}
$$
which is a smooth b-half density times a distribution conormal to the diagonal $\sigma = 1, y = y'$ uniformly to the boundary. 

In a neighbourhood of  the interior of $\bfo$ we can use coordinates 
\begin{equation}
\kappa = \frac{k}{x}, \ \kappa' = \frac{k}{x'}, \  y, \ y'  \text{ and } k.
\label{bfocoords}\end{equation}
 In these coordinates the kernel of the identity is 
$$
\delta(\kappa - \kappa') \delta(y-y') \left|\frac{d\kappa d\kappa' dy dy' dk}{k}\right|^{\demi}.
$$

Next consider the region near 
 $\bfo \cap \sca$ and in a neighbourhood of $\Delta_{k, \sca}$. The boundary hypersurface $\sca$ is obtained from $M^2_{k,b}$ by blowing up $\{ x/k = 0, x'/x - 1 = 0, y = y' \}$. The new coordinates are 
$$
X = \frac{x/x' - 1}{x/k} = k (\frac1{x'} - \frac1{x}), \ Y = k \frac{y-y'}{x}, \ \lambda = x/k, \ y, \ k. 
$$
In terms of these the kernel of the identity is
$$
\lambda^{-n/2} \delta(X) \delta(Y) \big| \frac{dX dY d\lambda dy dk}{\lambda k} \big|^{\demi}.
$$
Notice that $\lambda$ is a boundary defining function for $\sca$ locally near $\Delta_{k, \sca}$. 
We see that the kernel of the identity on half-densities is $\rho_{\sca}^{-n/2}$ times $\nu$ times a distribution conormal to $\Delta_{k, \sca}$ uniformly to the boundary of $\MMksc$.

Next consider the kernel of the operator $P + k^2$. We let this operator act on half-densities according to the formula
$$
(P+k^2) ( f |dg|^{1/2} ) = ((P + k^2) f) |dg|^{1/2}.
$$
Another way of viewing this prescription is that  derivatives are endowed with the flat connection on the half-density bundle that annihilates the scattering half-density $|dg|^{1/2}$. However, we wish to write the kernel of $P + k^2$ with respect to $\nu = |dg_b dg_b' dk/k|^{1/2}$. Let us therefore use the flat connection on half-densities that annihilates the b-half density $|dg_b|^{1/2}$. In terms of this, the kernel of $P + k^2$ is (in coordinates $(x,y)$)
\begin{equation}\begin{gathered}
 x^{-n/2} \Big( -(x^2 \pl_x)^2  - (n-1) x^3 \partial_x + x^2 \Delta_{S^{n-1}} + \tilde W + k^2 \Big) x^{n/2} \otimes \nu \\ = (x \, P_b \, x)  \otimes \nu.
\end{gathered}\label{PPb}\end{equation}
Notice that the vector field $x^2 \partial_x$ acting on the left factor lifts to $\MMksc$ to be transverse to $\Delta_{k, \sca}$ except at $\bfo$ where it vanishes to first order. Equivalently, the b-vector field $x \partial_x$ is transverse to $\Delta_{k, \sca}$ except at $\sca$ where it has growth of order $1$.  It follows that the kernel of $P + k^2$ is $\rho_{\bfo}^2$ times $\nu$ times a distribution which is conormal to $\Delta_{k, \sca}$ uniformly to the boundary.

\subsubsection{Polyhomogeneous functions and index sets}\label{pfis}
Let $X$ be a manifold with corners and $\mc{F}$ its set of boundary hypersurfaces. Let us recall that the index family consists of a subset $\mathcal{E}_{\rm f}$ of $\CC \times \NN$ for each ${\rm f}$ in the set  $\mc{F}$ of boundary hypersurfaces of $M^2_{k, \sca}$, satisfying some conditions given below. 
We say that a function $f$ on a manifold with corners $X$ is polyhomogeneous conormal (polyhomogeneous or phg for short) at the boundary with index family $\mc{E}$ if it is conormal  (i.e. if it remains in a fixed weighted $L^2$ space under repeated application of vector fields tangent to the boundary of $X$) 
and if  for each $s \in \RR$ we have
$$
\Big( \prod_{H \in \mc{F}} \prod_{(z, p) \in \mc{E}_H \text{ s.t. } \Re z \leq s} (V_H - z) \Big) f = O\big( (\prod_{H} \rho_H )^s \big)
$$
where  $V_H$ is a smooth vector field on $X$ that takes the form $V_H = \rho_H \partial_{\rho_H} + O(\rho_H^2)$ near $H$. This implies that $f$ has an asymptotic expansion in powers and logarithms near each boundary hypersurface. In particular, near the interior of $H$, we have 
$$
f = \sum_{{z,p} \in \mc{E}_H \text{ s.t. } \Re z \leq s} a_{(z,p)} \rho_H^z (\log \rho_H)^p +O(\rho_H^s)
$$
for every $s \in \RR$, where $a_{(z,p)}$ is smooth in the interior of $H$ (in fact, $a_{(z,p)}$ will itself be polyhomogeneous on $H$). 

To make sense of this definition, index  sets $\mc{E}_{\rm f}$ are required to have the properties that for each $M \in \RR$ the number of points $(\beta, j) \in \mc{E}_{\rm f}$ with $\Re \beta \leq M$ is finite, that if $(\beta, j) \in \mc{E}_{\rm f}$ then $(\beta + 1, j)\in \mc{E}_{\rm f} $ and that if $j > 0$ then also $(\beta, j-1) \in \mc{E}_{\rm f}$. These conditions ensure that the space of phg functions is well-defined (independent of the choice of $V_H$)  and closed under multiplication by $C^\infty(X)$. We remark that the index family consisting of the index set $\NN \times \{ 0 \}$ at each boundary hypersurface corresponds precisely to smooth functions on $M^2_{k, \sca}$. 

Recall the operations of addition and extended union of two index sets $E_1$ and $E_2$, denoted $E_1 + E_2$ and $E_1 \extunion E_2$ respectively:
\begin{equation}\begin{gathered}
E_1 + E_2 = \{ (\beta_1 + \beta_2, j_1 + j_2) \mid (\beta_1, j_1) \in E_1 \text{ and } (\beta_2 , j_2) \in E_2 \} \\
E_1 \extunion E_2 = E_1 \cup E_2 \cup \{ (\beta, j) \mid \exists (\beta, j_1) \in E_1, (\beta, j_2) \in E_2 \text{ with } j = j_1 + j_2 + 1 \}.
\end{gathered}\end{equation}

For what follows, we write $q$ for the index set 
\begin{equation}
\{ (q + n, 0) \mid n = 0, 1, 2, \dots \}
\label{short}\end{equation} for any $q \in \RR$, and we   define $\mc{N}$ and $\mc{N}_2$ to be the index sets
\begin{equation}\label{defn}
\mc{N}:=\{(k,l)\in\nn_0\x\nn_0; l\leq k\}, \ \mc{N}_2:=\{(k,l)\in\nn_0\x\nn_0; l\leq (k+1)^2/4\}.
\end{equation} 
For any index set $E$ and $q \in \RR$, we write $E \geq q$ if $\Re \beta \geq q$ for all $(\beta, j) \in E$ and if $(\beta, j) \in E$ and $\Re \beta = q$ implies $j = 0$. We write $E > q$ if there exists $\epsilon > 0$ so that $E \geq q + \epsilon$. Finally we say that $E$ is integral if $(\beta, j) \in E$ implies that $\beta \in \mathbb{Z}$.

\subsection{Operators defined on $M^2_{k, \sca}$}
We construct a parametrix $G(k)$ for $R(k)$ within a space of operators $\Psi^{m, \mc{E}}_k(M; \Omegab)$ given by half-density kernels on $M^2_{k, \sca}$. To each such kernels we can canonically associate an operator depending parametrically on $k$,  as follows: we write the kernel as $K \otimes \nu $ where $K$ is a distribution and $\nu = |dg_b dg'_b dk/k|^{1/2}$; then $( K|_{\{ k = k_1 \}}) \otimes |dg_b dg'_b|^{1/2}$ defines an operator on half-densities on $M$ for every $k_1 \in (0, k_0)$. 
For example, the resolvent of the Laplacian on $\RR^3$ is represented by
\begin{equation}
\frac1{4\pi} \frac{e^{-k|z-z'|}}{|z-z'|} \otimes \big| dz dz'|^{1/2} \otimes\big|  \frac{dk}{k} \big|^{1/2} =  \frac{(x x')^{-3/2}}{4\pi} \frac{e^{-k|z-z'|}}{|z-z'|} \otimes \nu
\label{3d}\end{equation}

We want a `variable coefficient' class of operators modelled on the free resolvents on $\RR^n$. The essential properties we shall require of our kernels are (i) conormality at the diagonal, uniformly  up to the boundary, and (ii) polyhomogeneity at the boundary away from the diagonal. We remark that the diagonal $\Delta_{k, \sca} \subset M^2_{k, \sca}$ is a p-submanifold transverse to each boundary hypersurface so conormality makes perfect sense here. 

The right-hand expression in \eqref{3d} for the resolvent on $\RR^3$ shows that the kernel vanishes to order $-2$ at $\bfo$ and order $-n/2$ at $\sca$. Moreover, the kernel vanishes to infinite order at $\bfc$, $\lb$ and $\rb$. This motivates the following definitions:

\begin{defn}\label{MMksc-op}
The class $\Psi^{m, \mc{E}}_k(M; \Omegab)$ of operators on half-densities on $M$ is the set of half-density kernels $K$ on $\MMksc$ that can be written as a sum of two terms $K = K_1 + K_2$, where 

(i) $\rho_{\sca}^{n/2} K_1$ is  supported near and conormal to $\Delta_{k, \sca}$ of order $m$, uniformly up to the boundary, and 

(ii) $\rho_{\sca}^{n/2}K_2$ is smooth in the interior of $M^2_{k, \sca}$ and classical conormal at the boundary with respect to some index family $\mathcal{E}$, i.e. an index set $\mc{E}_{\textrm{f}}$  at each boundary hypersurface $\textrm{f}$. Moreover, we assume that
\begin{equation}
\mc{E}_{\bfc}, \ \mc{E}_{\rb} \text{ and } \mc{E}_{\lb} \text{ are empty }
\label{empty}\end{equation}
and that (in the notation of \eqref{short})
\begin{equation}
\mc{E}_{\zf}, \mc{E}_{\bfo} \text{ and }  \mc{E}_{\sca} \text{ contain the $C^\infty$ index set } 0.
\label{consistency}\end{equation}
\end{defn}

\begin{defn}
The class $\Psi^{m, (a_{\bfo}, a_{\zf}, a_{\sca}); \mc{E}}_k(M; \Omegab)$ of operators on half-densities on $M$ is the class of half-density kernels $K$ that can be written as a sum of two terms $K = \tilde K_1 + K_2$, where $\tilde K_1 = \rho_{\bfo}^{a_{\bfo}}\rho_{\zf}^{a_{\zf}}\rho_{\sca}^{a_{\sca}} K_1$ and $K_1$ and $K_2$ are in Definition~\ref{MMksc-op}, with \eqref{consistency} replaced by
\begin{equation}
\mc{E}_{\bfo}, \mc{E}_{\zf} \text{ and }  \mc{E}_{\sca} \text{ contain the  index set } a_{\bfo}, \ a_{\zf}  \text{ and } a_{\sca} \text{ respectively.}
\label{consistency-2}\end{equation}

\end{defn}

\begin{remark} The order $-n/2$ of growth of the kernel at $\sca$ reflects the fact that we are expressing the kernel in terms of b-half-densities when scattering half-densities would be more natural here. \end{remark}

\begin{prop}\label{prop:composition} Suppose that $$A \in  \Psi_k^{m, (a_{\bfo}, a_{\zf}, a_{\sca}); \mc{A}}(M; \Omegab) \text{ and } B \in  \Psi_k^{m', (b_{\bfo}, b_{\zf}, b_{\sca});\mc{B}}(M; \Omegab)$$ satisfy \eqref{empty} and \eqref{consistency-2}. Then the composition $A \circ B$ is an element of $$\Psi_k^{m+m', (a_{\bfo} + b_{\bfo}, a_{\zf} + b_{\zf}, a_{\sca} + b_{\sca});\mc{C}}(M; \Omegab)$$  with 
\begin{equation}\begin{aligned}
\mc{C}_{\sca} &= \mc{A}_{\sca} + \mc{B}_{\sca} \\
\mc{C}_{\zf} &= \big( \mc{A}_{\zf} + \mc{B}_{\zf} \big) \extunion \big(  \mc{A}_{\rbo} + \mc{B}_{\lbo} \big) \\
\mc{C}_{\bfo} &= \big( \mc{A}_{\lbo} + \mc{B}_{\rbo} \big) \extunion \big(  \mc{A}_{\bfo} + \mc{B}_{\bfo} \big) \\
\mc{C}_{\lbo} &= \big( \mc{A}_{\lbo} + \mc{B}_{\zf} \big) \extunion \big(  \mc{A}_{\bfo} + \mc{B}_{\lbo} \big) \\
\mc{C}_{\rbo} &= \big( \mc{A}_{\zf} + \mc{B}_{\rbo} \big) \extunion \big(  \mc{A}_{\rbo} + \mc{B}_{\bfo} \big) \\
\mc{C}_{\bfc} &= \mc{C}_{lb} = \mc{C}_{\rb} = \emptyset
\end{aligned}\label{comp-if}\end{equation}

\end{prop}

\begin{proof} See Section~\ref{comp-proof}.
\end{proof}

\begin{cor}\label{cor} Suppose $E \in \Psi_k^{m, (a_{\bfo}, a_{\zf}, a_{\sca}); \mc{E}}(M; \Omegab)$  where $m < 0$, each of $a_{\bfo}, a_{\zf}, a_{\sca}$ is positive, and $\mc{E}_{\textrm{f}} > 0$ for all $\textrm{f}$. Then for large enough $N$, $E^N$ is Hilbert-Schmidt for each $k > 0$ with Hilbert-Schmidt norm $\| E^N(k) \|_{HS} \to 0$ as $k \to 0$. In particular, $E^N$ is compact for $N$ large enough, uniformly (in the sense above) as $k \to 0$.
\end{cor}

\begin{remark}\label{iterates-away} The operator $E$ will arise as the error term in a parametrix construction, i.e. we will have
$$
(P + k^2) G(k) = \Id + E(k)
$$
where $G(k)$ and $E(k)$ are kernels on $\MMksc$. We shall say that $E = E(k)$ `iterates away' if it has the properties in the corollary above. The point is that in the Neumann series for $(\Id + E)^{-1}$, each term becomes sucessively milder, and the Neumann series can be summed asymptotically modulo a smooth kernel on $\MMksc$ that vanishes to infinite order at the boundary. Moreover a finite number of terms in the Neumann series gives an error that is compact uniformly in $k$.
\end{remark}

\begin{proof} Choose $\delta > 0$ so that $-m \geq \delta$, $a_{\bfo}, a_{\zf}, a_{\sca} \geq \delta$, and $\mc{E}_{\textrm{f}} \geq \delta$ for all $\textrm{f}$. Then, by Proposition~\ref{prop:composition}, $E^N$ has a similar property with respect to $N \delta$. For $N \delta > n$ the kernel is continuous and vanishes to order greater than $n$ at each boundary hypersurface of $\MMksc$ which implies that it is Hilbert-Schmidt for each $k > 0$ with Hilbert-Schmidt norm tending to zero  as $k \to 0$.
\end{proof}

\subsection{Normal operators}
Assume that 
\begin{equation}
\mc{E}_{\zf} , \mc{E}_{\bfo}  \text{ and } \mc{E}_{\sca} \geq 0. 
\label{nn}\end{equation}
Then the restriction of $A \in \Psi^{m, \mc{E}}(M)$ to each of these faces  is well-defined and can be viewed as an operator on functions on a model space. On $\zf$ this model operator is a b-pseudodifferential operator\footnote{This is not quite true, since there may be a polyhomogeneous expansion of $I_{\zf}(A)$ at the boundary of $\zf$ which is not allowed in the b-calculus as defined in \cite{APS}. However this is an inessential point.} acting on half-densities on $M$, on $\bfo$ it is a pseudodifferential operator acting on half-densities defined on $S^{n-1} \times (0, \infty)$ and on $\sca$ it is a family of convolution pseudodifferential operators, parametrized by $\hat z \in S^{n-1}$ and by $k \in (0, k_0]$, acting on half-densities defined on $\RR^n$. We call the model operator at $\textrm{f} \in \{ \zf, \bfo, \sca \}$ the normal operator of $A$ at $\textrm{f}$ and denote it $I_{\textrm{f}}(A)$.

\begin{prop}
The normal operators respect  composition: if $A \in \Psi^{m, \mc{A}}(M, \Omegab)$ and $B \in \Psi^{m, \mc{B}}(M, \Omegab)$ satisfy \eqref{nn} and 
\begin{equation}
\mc{A}_{\rbo} + \mc{B}_{\lbo} > 0 \text{ and } \mc{A}_{\lbo} + \mc{B}_{\rbo} > 0
\label{pos-cond}\end{equation} then
\begin{equation}\begin{aligned}
I_{\zf}(A \circ B) &= I_{\zf}(A) \circ I_{\zf}(B) \\
I_{\bfo}(A \circ B) &= I_{\bfo}(A) \circ I_{\bfo}(B) \\
I_{\sca}(A \circ B) &= I_{\sca}(A) \circ I_{\sca}(B) 
\end{aligned}\end{equation}
\end{prop}

\begin{proof} The third composition property is just that of the scattering calculus with a smooth parameter $k$. 

To prove the first identity, note that $I_{\zf}(A \circ B)$ is the pushforward of $\pi_L^* A \otimes \pi_R^* B \otimes \pi_C^* \nu \otimes | dk/k|^{-1/2}$ restricted to the boundary hypersurface of $\MMMksc$ which is the lift of $M^3 \times \{ 0 \}$. This boundary hypersurface is canonically isomorphic to $M^3_b$ and the map to $\zf$ is the canonical b-fibration. This implements the product in the b-calculus, so the first line follows immediately. (We remark that there is another boundary hypersurface of $\MMMksc$ that projects under $\pi_C$ to $\zf$ but under the assumption \eqref{pos-cond} it does not contribute to the leading order term.) 

The second identity follows similarly. 
\end{proof}

\begin{remark} Following on from the last point in Remark~\ref{essential}, the face $\bfo$ is the space $(\partial M \times I)^2$, where $I$ is the compactification of the interval $(0, \infty)_\kappa$ parametrized by $\kappa = k/x$, with a `b'-blowup at $\kappa = \kappa' = 0$ and a scattering blowup at $\kappa = \kappa' = \infty$. Looking ahead to Section~\ref{nokernel}, we see from \eqref{bfo-normal-1} and \eqref{bfo-normal-2} that the normal operator of $P$ at $\bfo$ is of b-type at $\kappa = 0$ and of scattering type at $\kappa = \infty$, so we see the transition between b- and scattering type behaviour quite explicitly. 
\end{remark}

\subsection{Strategy}\label{strategy}

We attempt to solve the equation
\begin{equation}\label{tosolve}
(P+k^2) G(k) = \kappa_{\rm Id}
\end{equation}
where $G = G(k)$ is in $\rho_{\bfo}^{-2} \Psi^{-2, \mc{E}}(M^2_{k, \sca}; \Omegab)$. Of course we cannot do this in one step, so we first try to write down a parametrix $G$ with an error term $E = E(k)$ which iterates away in the sense of Remark~\ref{iterates-away}. 

In the simplest case that $P$ has no $L^2$ null space nor a zero-resonance, our parametrix $G$ will have the index set $0$ at $\zf$ and the index set $-2$ at $\bfo$. Let us write $G_{\zf}^0$ for the restriction of $G$ to $\zf$ (which is canonically defined since $G$ is a b-half-density) and $G_{\bfo}^{-2}$ for the restriction of $k^2 G$ to $\bfo$. In general, if we need to consider more terms in the Taylor series of $G$  we shall write $G_{\textrm{f}}^j$, $j \in \NN$, for the $j$th term in the Taylor series of $G$ at any face $\textrm{f}$. Of course we need to specify a boundary defining function for $\textrm{f}$ if $j = 0$. \emph{We shall always use the boundary defining function $k$ at $\zf$, $\rbo$, $\bfo$ or $\lbo$}; this is the most convenient choice since it commutes with the operator $P + k^2$. We remark that $k$ is only a boundary defining function in the \emph{interior} of each of these faces, but this is of no importance; we only have to remember that the coefficients $G_{\textrm{f}}^j$ then may become more and more singular at the boundary of $\textrm{f}$ as $j$ increases. 

With this notation, then we see that if $G \in \rho_{\bfo}^{-2} \Psi^{-2, \mc{E}}(M^2_{k, \sca}; \Omegab)$ then our requirements are that $E^0_{\zf}$, $E^0_{\bfo}$ and $E^0_{\sca}$ all vanish. This amounts to inverting model elliptic operators at each of these three boundary hypersurfaces. It turns out that in high dimensions, $n \geq 5$, this is sufficient to obtain an error term that iterates away. In low dimensions or in the presence of $L^2$ null space or a zero-resonance we need to work harder, that is, to solve further model problems. We shall begin with the simplest case first.


\section{Resolvent kernel  when
$(M,g)$ is asymptotically Euclidean}\label{nokernel}

Let $M,g)$ be an asymptotically Euclidean manifold. We shall construct a para\-metrix for $(P + k^2)^{-1}$ with an error that `iterates away' in the sense of Remark~\ref{iterates-away}. Throughout this paper we assume that $P$ has neither nontrivial $L^2$ null space, nor a zero-resonance.  Initially we shall assume that $n \geq 5$, as the parametrix construction then is as straightforward as possible. In section~\ref{34} we treat the cases $n = 3,4$.

\subsection{Singularity at the diagonal $\Delta_{k, \sca}$.}\label{singdiagonal} The operator $P + k^2$ is elliptic as an element of $\Psi^{2, (2, 0, 0); *}(M, \Omegab)$, i.e. its symbol multiplied by $\rho_{\bfo}^{-2}$ is elliptic uniformly on $N^* (\Delta_{k, \sca})$. We specify the symbol of $G(k)$ to be the inverse (in the sense of the noncommutative product corresponding to operator composition) of that of $P + k^2$ modulo symbols of order $-\infty$. This is consistent with membership of $G(k)$ in $\Psi^{-2; (-2, 0, 0), *}(M, \Omegab)$ and determines the symbol of $G(k)$ modulo symbols of order $-\infty$. 

\subsection{Term at $\sca$.} The normal operator of $P + k^2$ at $\sca$ is the family of operators $\Delta + k^2$ acting on half-densities on $\RR^n$, parametrized by $\hat z \in S^{n-1}$ and $k \in (0, k_0]$. We specify the normal operator $I_{\sca}(G(k))$ to be the inverse operator $(\Delta + k^2)^{-1}$ which is well-defined for $k > 0$. 

\subsection{Term at $\zf$}\label{zfnokernel} 
Let us start with the face $\zf$, which is canonically identified with the b-blow-up $M^2_b:=[\Mbar\x\Mbar;\pl\Mbar\x\pl\Mbar]$.
On this face $\zf$, the equation is really a b-elliptic equation in disguise: we see from \eqref{PPb} that $P = x P_b x$ (where we interpret $P$ and $P_b$ as acting on half-densities, and $P$ annihilates the scattering half-density while $P_b$ annihilates the b-half-density)  where 
the operator $P_b$ is self-adjoint with respect to the b-metric $g_b$ and elliptic in the b-calculus. 
Moreover, its normal operator $I_{\rm ff}(P_b)$ is invertible. Therefore,  by the theory of b-elliptic operators given by Melrose \cite[Sec. 5.26]{APS}, there is a generalized inverse $Q_b$, a b-pseudodifferential operator of order $-2$ for the operator $P_b$ on $L^2_b$, such that 
\begin{equation}\label{pbqb}
P_b Q_b = Q_b P_b = \Id -\Pi_b.
\end{equation}
where  $\Pi_b$ is orthogonal projection on the $L_b^2$ kernel of $P_b$.
Fr{o}m \cite{APS}, an element $\varphi$ in the range of $\Pi_b$ is in $x^{\ndemi-1}L^\infty(\Mbar,\Omega^\demi_b)$ (see Theorem~\ref{rit} and Proposition~\ref{prop:complementary}), for all $\epsilon > 0$. This implies that $x^{-1}\varphi\in x^{n-2}L^2(\Mbar,\Omega^{\demi}_{\sca})$ satisfies $P (x^{-1}\varphi)=0$. So $x^{-1} \varphi$ is either a zero-mode or a zero-resonance. By our assumptions on $P$, we have 
\[\Pi_b=0.\]
Again from \cite{APS}, the kernel $Q_b$ is  conormal at the b-diagonal $\Delta_{k, \sca} \cap\zf$, uniformly up to the front face $\zf \cap \bfo$, and is the sum of three terms: one is supported close to and conormal at the b-diagonal,  uniformly to the boundary;  one is polyhomogeneous on $M^2_b$ with the index set  
$\mc{E}(Q_b)=(\mc{E}_{\rm ff}(Q_b), \mc{E}_{\rb}(Q_b), \mc{E}_{\lb}(Q_b))$ where $\rb,\lb,{\rm ff}$ denote
the right boundary, left boundary and front face of the blow-up $M^2_b$, and 
\[\mc{E}_{\rb}(Q_b)=\mc{E}_{\lb}(Q_b)=\{ (n/2-1+l, k); 0\leq k\leq l\in\nn_0\}=\ndemi-1+\mc{N},\quad 
\mc{E}_{\rm ff}=0\]  
and the other term is the lift of a polyhomogeneous distribution on $\Mbar\x\Mbar$, with index sets $\mc{E}_{\rb}(Q_b)$ at the right and left boundaries.

So we have 
$$
P_b Q_b = \Id.
$$

Therefore we set
$$
G_{\zf}^0 = (x x')^{-1} Q_b.
$$

\subsection{Leading term at $\bfo$}\label{bfonokernel}
We use coordinates (\ref{bfocoords})
near the interior of $\bfo$ in terms of which we have $\textrm{int}(\bfo)=\{k=0,\kappa,\kappa'\in (0,\infty)\}$.
The operator $P+k^2$ in these coordinates reads
\[P+k^2=k^2\kappa^{-\ndemi-1} \Big(-(\kappa\pl_\kappa)^2+\Delta_{S^{n-1}}+(n-2)^2/4+\kappa^2+W
\Big)\kappa^{\ndemi-1}.\]
Since $W\in x\textrm{Diff}_b(M)$, the first germ of $P+k^2$ at $\bfo$ is at order $2$, we consider
\[I_{\bfo}(k^{-2}(P+k^2))=\kappa^{-\ndemi-1}\Big(-(\kappa\pl_\kappa)^2+\Delta_{S^{n-1}}+(n-2)^2/4+\kappa^2\Big)\kappa^{\ndemi-1}.\]  
Here we wrote $P$ with respect to the flat connection on half-densities annihilating $|dg|^{1/2}$. In terms of the  b-flat connection annihilating the b-half-density $|dg_b|^{1/2}$, our operator is 
\begin{equation}\label{bfo-normal-1}
\kappa^{-1}\Big(-(\kappa\pl_\kappa)^2+\Delta_{S^{n-1}}+(n-2)^2/4+\kappa^2\Big)\kappa^{-1} := \kappa^{-1} P_{\bfo} \kappa^{-1}.
\end{equation}
Setting $M_0:=(0,\infty)\x S^{n-1}$, we can view $\textrm{int}(\bfo)$ as $M_0\x M_0$. 
The operator $P_{\bfo}$ acting on the left on $M_0\x M_0$ on b half-densities
$f(\kappa,y)|\kappa^{-1}d\kappa dy|^{\demi}$ has an inverse 
in terms of spherical harmonics $(\phi_j(y))_j$ from Sturm-Liouville theory: 
if $H$ be the Heaviside function, $I_\nu,K_\nu$ the modified Bessel functions, then
\begin{equation}\label{Qbfo}
\begin{split}
 Q_{\bfo}:=&\sum_{j=0}^\infty\Pi_{E_j}\Big(I_{j+\ndemi-1}(\kappa)K_{\ndemi-1+j}(\kappa')H(\kappa'-\kappa)\quad\quad\quad\quad\\
&\quad\quad\quad\quad\quad+I_{j+\ndemi-1}(\kappa')K_{\ndemi-1+j}(\kappa)H(\kappa-\kappa')\Big) \left|\frac{d\kappa dyd\kappa'dy'}{\kappa\kappa'}\right|^\demi
\end{split}\end{equation}
solves on b-half densities on $M_0\x M_0$ 
$$
P_{\bfo} Q_{\bfo} = \Id.
$$

Therefore we set
\begin{equation}
G_{\bfo}^{-2} = (\kappa \kappa') Q_{\bfo} .
\label{Gbfo}\end{equation}
 
\subsection{Consistency} 
Before proceeding further we need to check that our specified values for $\sigma(G(k))$, $G^0_{\sca}$, $G^0_{\zf}$ and $G^{-2}_{\bfo}$ are consistent. 

We first check consistency between $G^0_{\zf}$ and $G^{-2}_{\bfo}$. 
To do this we must verify that the restrictions of $\rho_{\bfo}^{2} G_{\zf}^0$ and $\rho_{\zf}^2 G_{\bfo}^{2}$ agree at $\zf \cap \bfo$, where $\rho_{\zf} \cdot \rho_{\bfo} = k$. It is most convenient to take
$\rho_{\bfo}^{-2} = (x x')^{-1}$ and $\rho_{\zf}^2 = \kappa \kappa'$. Thus, we compare the restrictions of $(x x') G^0_{\zf}$ and $(\kappa \kappa')^{-1} G^{-2}_{\bfo}$ at the intersection $\bfo \cap \zf$. Fr{o}m $G^0_{\zf}$ we get the normal operator (in the sense of the b-calculus) of $Q_b$. This is the inverse of the normal operator of $P_b$  which has the exact expression $(sD_s)^2 + (n/2-1)^{2} + \Delta_{S^{n-1}}$ where $s:=\beta^*(x/x')\in\rr$ is a coordinate on the fibers of ff. The solution is given by standard Sturm-Liouville theory in terms of the eigendecomposition of $\Delta_{S^{n-1}}$ by 
\begin{equation}\label{nkernelff}
\sum_{j=0}^{\infty} \Pi_{E_j} \frac{e^{-(j+\ndemi-1) |\log s|}}{2j+n-2}|dydy'ds/s|^{\demi}.
\end{equation}

On the other hand, to compute the restriction of  $(\kappa \kappa') G^{-2}_{\bfo}$ at the intersection $\bfo \cap \zf$, we recall the definitions of $I_\nu,K_\nu$ and their asymptotics
\begin{equation}\label{inuknu1}
I_\nu(z)=2^{-\nu}z^{\nu}\sum_{k=0}^{\infty}\frac{4^{-k}z^{2k}}{k!\Gamma(\nu+k+1)}, \quad
K_{\nu}(z)=-\frac{\pi}{2} \frac{I_{\nu}(z)-I_{-\nu}(z)}{\sin \nu \pi}, \nu \notin \mathbb{Z},
\end{equation}
while for  $\nu = n \in \mathbb{Z}$,   $I_n$ is entire, $I_n = I_{-n}$ and 
$K_\nu$ is given by the limiting value of the right hand side  as $\nu \to n$. As $z \to 0$, we have
\begin{equation}\label{inuknu2}
I_\nu(z) \sim \frac1{\Gamma(\nu + 1)} (\frac{z}{2})^\nu, \quad
K_\nu(z) \sim \frac{\Gamma(\nu)}{2} (\frac{z}{2})^{-\nu},  \quad \nu > 0. 
\end{equation}
Also  $K_{\nu}(z)=O(e^{-c|z|})$ for some $c>0$ as $|z|\to \infty$ with $z,\nu\in\rr$. At $\bfo \cap \zf$ we have $\kappa = \kappa' = 0$ and only the leading behaviour \eqref{inuknu2} of $I_\nu$ and $K_\nu$ at $\kappa=0$ contributes to this restriction. We see from this that the restrictions agree identically. 

We next check consistency between $G^0_{\sca}$ and $G^{-2}_{\bfo}$. The face $\sca$ is fibred over $(k, y)$ with fibres having a Euclidean structure. Let $z$ be a Euclidean variable on the fibres of the scattering double space $M^2_{\sca}$; then $kz$ is a Euclidean variable on the fibres of $\sca$. 
We have defined $G^0_{\sca} = k^{n-2} f(kz) |dz|$, where $f(z)$ is the kernel of $(\Delta + 1)^{-1}$ on $\RR^n$ and hence $k^{n-2} f(kz)$ is the kernel of $(\Delta + k^2)^{-1}$ on $\RR^n$. 

At $\bfo$, let us write the normal operator $G^{-2}_{\bfo}$ in terms of the scattering half-density $k^n| dg dg' dk/k|^{1/2}$ which restricts to the half-density $| \kappa^{n-1} d\kappa dy {\kappa'}^{n-1} d\kappa' dy'|^{1/2}$ at $\bfo$. That is, let
$$
G^{-2}_{\bfo} = \tilde Q_{\bfo} k^n| dg dg' dk/k|^{1/2}.
$$
Correspondingly we use the flat connection on half-densities that annihilates the scattering half-density.  The PDE for $\tilde Q_{\bfo}$ then takes the form 
\begin{equation}\label{bfo-normal-2}
\Big( -\partial_{\kappa}^2 - \frac{n-1}{\kappa} \partial_{\kappa} + \kappa^{-2} \Delta_{S^{n-1}} + 1  \Big)  \tilde Q_{\bfo} = \kappa^{1-n} \delta(\frac1{\kappa} - \frac1{\kappa'}) \delta_y(y').
\end{equation}
This is just the equation $\Delta + 1$  on $M_0$ viewed as $\RR^n$ in polar coordinates $(\kappa, y)$, with $\kappa$ playing the role of the radial variable. It follows that 
we have $\tilde Q_{\bfo}  =f(\kappa y - \kappa' y')$, where we view $\kappa y$ and $\kappa' y'$ as lying in $\RR^n$. 

The boundary hypersurface $\bfo$ may be obtained from $M_0$ by 
first compactifiying $M_0$, then forming the b-stretched product, and then blowing up the boundary of the b-diagonal at the `$\infty$' end only. This last blowup is at $\sca \cap \bfo$ and this face is fibred over $S^{n-1}$ with fibres that are Euclidean with Euclidean variable $\tilde z = \kappa y - \kappa' y'$. It follows that the restriction of $\tilde Q_{\bfo}$ to this face is $f(\tilde z)$. We finally observe that $\tilde z$ agrees with the limiting value of the coordinate $kz$ taken from $\sca$.  This shows  the compatibility between $G^0_{\sca}$ and $G^{-2}_{\bfo}$.

Next we check consistency between the symbol of $G(k)$ and the normal operators. This is immediate because each of the normal operators at $\sca, \zf$ and $\bfo$ solve an elliptic pseudodifferential equation  so the symbol at the diagonal on each face is determined uniquely modulo symbols of order $-\infty$, and hence necessarily agree modulo symbols of order $-\infty$ with that coming from the interior of $\Delta_{k, \sca}$. 

We next consider the behaviour of these normal operators at the remaining boundary faces. We note that $G^0_{\sca}$ decays to infinity order at the boundary $\sca \cap \bfc$, uniformly down to $k=0$. (This is because a boundary defining function for $\bfc$ near $\sca \cap \bfc$ is $k|z-z'|$.) Similarly, the exponential decrease of $K_\nu(\kappa)$ as $\kappa \to \infty$ shows that $G^{-2}_{\bfo}$ decays rapidly at the intersection $\bfo \cap \bfc$, $\bfo \cap \lb$ and $\bfo \cap \rb$. Finally, both $G^{-2}_{\bfo}$ and $G^0_{\zf}$ are polyhomogeneous at $\lbo$ and $\rbo$ with index set $n/2 - 2+\mc{N}$. Moreover $G_{\bfo}^{-2}$ has index set at $\zf$ given by
\begin{equation}\label{indexgbfo}
0 \text{ if }n \text{ is odd}, \quad 0+\{(n-2+2l,1) \mid l\in \nn_0 \} \text{ if }n\text{ is even} .\end{equation}

We conclude that there exists a $G = G(k) \in \Psi^{-2, (-2, 0, 0); \mc{G}}(M, \Omegab)$ satisfying $\mc{G}_{\sca} =0$, $\mc{G}_{\bfo} = -2$, $\mc{G}_{\lbo} = \mc{G}_{\rbo} = n/2 - 2+\mc{N}$, $\mc{G}_{\zf}$ is given by (\ref{indexgbfo}), the index set at $\rb,\lb,\bfa$ is empty, and $G$ has $G^{-2}_{\bfo},G^0_{\zf},G^{0}{\sca}$ as leading term
at $\bfo,\zf,\sca$, with symbol at the diagonal given in \ref{singdiagonal}. 

\subsection{Exact resolvent}
Applying $P + k^2$ to our parametrix $G(k)$, we obtain 
$$
(P + k^2) G(k) = \Id + E(k)
$$
where $E(k) \in \Psi^{-\infty, (1, 1, 1); \mc{E}}(M, \Omegab)$, where $\mc{E}$ satisfies $\mc{E}_{\zf} \geq 1$, $\mc{E}_{\bfo} \geq 1$, $\mc{E}_{\sca} \geq 1$, $\mc{E}_{\lbo} \geq n/2$, $\mc{E}_{\rbo} \geq n/2 - 2$ and \eqref{empty}. (We emphasize that since the operator $P + k^2$ vanishes to order $2$ at $\bfo$, we \emph{automatically} gain two orders at $\bfo$, and the fact that we solved the model problem at $\bfo$ gives us the additional order of vanishing. Hence $E(k)$ is three orders better at $\bfo$ than $G(k)$.) 
In particular, 
$\mc{E}_{\textrm{f}} > 0$ for each boundary hypersurface $f$. It follows from Corollary~\ref{cor} that $E(k)$ is Hilbert-Schmidt with Hilbert-Schmidt norm going to zero as $k \to 0$. It follows that $\Id + E(k)$ is invertible for small $k$. Moreover, $E(k)$ `iterates away' in the sense of Remark~\ref{iterates-away}.  It follows from standard arguments that $(\Id + E(k))^{-1}$ lies in the calculus. 
Inspection of \eqref{comp-if} and an inductive argument shows that for $l\geq 2$ 
\[ E(k)^l=\in \Psi_k^{-\infty,(l,l,l);\mc{E}^l}(M,\Omegab)\]
where $\mc{E}^l$ is the index set
\[ \mc{E}_{\sca}^l=l, \quad \mc{E}_{\zf}^l, \mc{E}_{\bfo}^{l}\subset l+\mc{N}_l, \]
\[
\mc{E}_{\rbo}^l\subset \ndemi+l-3+\mc{N}_l,\quad  \mc{E}^l_{\lbo}\subset \ndemi+l-1+\mc{N}_l\] 
with $\mc{N}_l:=\{(k,j)\in\nn_0\x\nn_0; j\leq l(k+1)\}$.
Let us write 
$$
(\Id + E(k))^{-1} = \Id + S(k).
$$
Then $S(k)$  is an element of $\Psi^{-\infty, (1,1,1); \mc{E}'}(M, \Omegab)$ for some index set $\mc{S}$ determined by $\mc{E}$. Inspection of \eqref{comp-if} and an inductive argument shows that 
\[\mc{S}_{\zf}, \mc{S}_{\bfo}\subset 1\cup (2+\mc{N}^2),\quad 
\mc{S}_{\lbo} \subset (\frac{n}{2} + \mc{N}^2),  \quad
\mc{S}_{\rbo}\subset (\ndemi-2) \cup (\ndemi-1+\mc{N}^2),\quad\mc{S}_{\sca} = 1 \]
where $\mc{N}^2$ is defined by
\[\mc{N}^2:= \{ (k,j)\in\nn_0\x\nn_0; j\leq (k+2)^2/4 \}.\]

It follows that the resolvent kernel can be written 
$R(k) = G(k) + G(k) S(k)$ and also lies in the calculus. A final application of the pushforward theorem gives \eqref{res-kernel}, with 
\begin{equation}
\begin{gathered}
\mc{R}_{\zf} \subset 0 \cup (1+\mc{N}^2),\quad
\mc{R}_{\lbo},\mc{R}_{\rbo}  \subset (\frac{n}{2} - 2+ \mc{N}^2), \\
\mc{R}_{\bfo}\subset -2\cup (-1+\mc{N}^2), \quad
\mc{R}_{\sca} = 0.
\end{gathered}\label{rif}\end{equation}
Moreover the correction term $G(k) S(k)$ vanishes to first order at $\zf$ and $\bfo$ so our parametrix is exact to leading order at these faces. 
Moreover,  the asymptotic expansion of the resolvent $R(k)$ as $k \to 0$ can, in principle, be generated to any given order at each boundary hypersurface from computing a finite number of terms in the Neumann series, i.e. from $G(k) (\Id - E(k) + E(k)^2 - \dots + (-1)^N E(k)^N)$ for sufficiently large $N$. 

\begin{remark} We emphasize that the index families given above are not sharp as far as the logarithmic exponents are concerned. 
\end{remark}

\subsection{Leading term at $\rbo$ and $\lbo$.}\label{34} We can improve the accuracy of our para\-metrix by specifying  $G^{n/2-2}_{\rbo}$ and $G^{n/2-2}_{\lbo}$ appropriately. When this is done correctly we obtain a better error term $E(k)$ with leading behaviour at $\rbo$ and $\lbo$ at order $n/2 - 1$ instead of $n/2 - 2$. This is useful for two reasons. The first is that this allows us to construct the resolvent in dimensions 3 and 4. 
The second is that it allows us to determine the leading behaviour of the resolvent itself at $\rbo$ and $\lbo$, which turns out to be crucial in understanding the $L^p$ boundedness properties of the Riesz transform. 

The term $G^{n/2-2}_{\rbo}$ must solve the equation $P_b ( x G^{n/2-2}_{\rbo}) = 0$ and match correctly with $G^0_{\zf}$ and $G^{-2}_{\bfo}$. Let us first consider the behaviour of $G^0_{\zf} = (x x')^{-1} Q_b$ near $\zf \cap \rbo$. 
Localizing near $\rbo$, we have by Theorem~\ref{reg}
\begin{equation}\label{qbas}
Q_b=\Big(v(z){x'}^{\ndemi-1}+O({x'}^{\ndemi})\Big)\left|\frac{dxdx'dydy'}{xx'}\right|^{\demi}\end{equation}
for some $v(z)$ with a polyhomogeneous expansion in $z = (x,y)$ as $x \to 0$. In order to match with the normal operator at the front face (here $\zf \cap \bfo$) $v$ can grow no faster than $x^{-n/2 + 1}$ as $x \to 0$. Proposition~\ref{prop:complementary} shows that there is a unique element of the null space of $P_b$ (up to scaling) with this property, which has the behaviour $c x^{-n/2 + 1} + O(x^{-n/2 + 2})$ as $x \to 0$ where $c$ is constant in $y$. Matching with the normal operator of $Q_b$ then gives the leading asymptotic 
\[v(z)=\frac{1}{(n-2){\rm Vol}(S^{n-1})}x^{-\ndemi+1}+O(x^{-\ndemi+2}).\]
To match with $G^0_{\zf}$ then we take
$$
G^{n/2-2}_{\rbo} = x^{-1} v(z) f(\kappa', y') \left|\frac{dxdx'dydy'}{xx'}\right|^{\demi}
$$
where $f(\kappa', y) \sim (\Gamma(n/2-1)2^{n/2-2})^{-1}{\kappa'}^{-n/2 - 2}$ as $\kappa' \to 0$ and otherwise is determined by the matching condition with $G^{-2}_{\bfo}$. At the intersection $\rbo \cap \bfo$, $\kappa = 0$ and $\kappa' \in (0, \infty)$. In particular we have $\kappa < \kappa'$. Examining 
\eqref{Qbfo} and \eqref{Gbfo} we see that the leading behaviour of $G^{-2}_{\bfo}$ at $\bfo \cap \rbo$ is 
$$
\frac1{\Gamma(n/2 - 1)} 2^{-n/2 +2} \kappa^{n/2} K_{n/2 - 1}(\kappa').
$$
Hence we set $f(\kappa', y)  =(\Gamma(n/2 - 1))^{-1} 2^{-n/2 +2}  \kappa' K_{n/2 - 1}(\kappa')$, i.e. 
$$
G^{n/2-2}_{\rbo} = \frac{ x^{-1} v(z) \kappa' K_{n/2 - 1}(\kappa')}{\Gamma(n/2 - 1) 2^{n/2 -2}}   \left|\frac{dxdx'dydy'}{xx'}\right|^{\demi}
$$
 and this matches with both  $G^0_{\zf}$ and $G^{-2}_{\bfo}$.
 
 We now set 
 $$
 G^{n/2-2}_{\lbo} = \frac{ {x'}^{-1} v(z') \kappa K_{n/2 - 1}(\kappa)}{\Gamma(n/2 - 1) 2^{n/2 -2}}   \left|\frac{dxdx'dydy'}{xx'}\right|^{\demi}
$$
to satisfy formal self-adjointness at this order. We find that the error term now is
$E(k) \in \Psi^{-\infty, (1, 1, 1); \mc{E'}}(M, \Omegab)$, where $\mc{E}'$ satisfies for some $\mc{E}'_{\zf} \geq 1-\eps$, $\mc{E}'_{\bfo} \geq 1-\eps$, $\mc{E}'_{\sca} \geq 1-\eps$, $\mc{E}'_{\lbo} \geq n/2 + 1-\eps$, $\mc{E}'_{\rbo} \geq n/2 - 1-\eps$ and \eqref{empty} for some small $\eps<1$. Then one can check that $G(k) (E(k))^N$ has an index set at $\rbo$ which is always $\geq n/2 - 1-\eps$. That means that the resolvent kernel itself has an index set at $\rbo,\lbo$ included in $(n/2-2)\cup \mc{Z}$ for some index set $\mc{Z}\geq n/2-1-\eps$, and the term of order $n/2-2$
satisfies $R^{n/2-2}_{\rbo} = G^{n/2-2}_{\rbo}$, i.e. we have found the actual leading behaviour of the resolvent kernel at $\rbo$. Moreover, the construction now goes through in dimensions 3 and 4, assuming that the $L^2$ null space is trivial and that there are no zero-resonances, since now the error term $E(k)$ vanishes to a positive order $n/2 - 1$ at the right boundary as opposed to $n/2 - 2$ as previously.\\ 

\begin{remark}\label{vzconstant}
We recall that $v(z)$ solves $Px^{\ndemi-1}v=0$ and $x^{\ndemi-1}v(z)$ is asymptotically constant at infinity. 
If $M$ has only one end, then Theorem~\ref{rit} implies that this is the only function (up to scaling) satisfying these properties. Hence, if $M$ has only one end and $V \equiv 0$, then $x^{\ndemi-1}v(z)$ is necessarily
constant. 
\end{remark}

\section{Resolvent kernel when $(M,g)$ is asymptotically conic}\label{nokernelconic}
 
We essentially give the differences with the asymptotically Euclidean case.
First recall $N_{\pl}:=\{(n-2)/2=\nu_0<\nu_1\leq\dots  \}$ is defined in (\ref{npl}).

\subsection{Singularity at the diagonal $\Delta_{k, \sca}$.}  This is completely similar to the asymptotically Euclidean case.

\subsection{Term at $\sca$.} The fibers of $\sca$ are Euclidean spaces radially compactified
and the normal operator of $P + k^2$ at $\sca$ is the family of flat Laplacians $\Delta + k^2$ acting on half-densities on these Euclidean spaces, parametrized by $y' \in \pl\Mbar$ and $k \in (0, k_0]$. 
Indeed if $y_1,\dots,y_n$ are local coordinates on $\pl\Mbar$, the functions $z_0=k/x-k/x',z_i:=k(y_i/x-y'_i/x')$ 
give Euclidean coordinates on the fibers $\sca$ and the normal operator of $\Delta_{g}+k^2$ at the fiber 
of $\sca$ with basis point $(k,y')$
at $\sca$ is $k^2(\Delta_z+1)$ where $\Delta_z$ is the Euclidean Laplacian $-\pl^2_{z_0}-\sum_{i,j}h^{ij}_0(y')\pl_{z_i}\pl_{z_j}$, $h_0(y')$ being the Euclidean metric 
on this fiber induced by $h_0$.
We specify the normal operator $I_{\sca}(G(k))$ to be, in each fiber, the inverse operator $k^{-2}(\Delta_z + 1)^{-1}$ which is well-defined for $k > 0$. 

\subsection{Term at $\zf$} 
As in the case $\pl\Mbar=S^{n-1}$, the theory of b-elliptic operators given by Melrose \cite[Sec. 5.26]{APS} shows that there is a generalized inverse $Q_b$, a b-pseudodifferential operator of order $-2$ for the operator $P_b$ on $L^2_b$, such that 
\begin{equation*}
P_b Q_b = Q_b P_b = \Id -\Pi_b.
\end{equation*}
where  $\Pi_b$ is orthogonal projection on the $L_b^2$ kernel of $P_b$. 
As in Section~\ref{nokernel}, due to our assumptions on $P$,  we have 
\[\Pi_b=0.\]
As in Section~\ref{zfnokernel},
$Q_b$ is  conormal at the b-diagonal $\Delta_{k, \sca} \cap\zf$, uniformly up to $\zf \cap \bfo$, and is the sum of three terms as described there, except that the polyhomogeneous term has the index set  
$\mc{E}(Q_b)=(\mc{E}_{\rm ff}(Q_b), \mc{E}_{\rb}(Q_b), \mc{E}_{\lb}(Q_b))$ where 
\[\mc{E}_{\rb}(Q_b)=\mc{E}_{\lb}(Q_b)\subset\{ (\nu, k)\in \rr\x\nn_0;  \nu\in N, k+1\leq \sharp\{\mu\in N_\pl;\mu\leq \nu\}\},\]
\[ 
\mc{E}_{\rm ff}=\nn_0\x \emptyset.\]  
So we have $P_b Q_b = \Id$ and we set
$G_{\zf}^0 = (x x')^{-1} Q_b.$

\subsection{Leading term at $\bfo$}
This is essentially similar to the case $\pl\bar{X}=S^{n-1}$ except that we have to replace $\Delta_{S^{n-1}}$ by $\Delta_{\pl\Mbar}$ and $P_{\bfo}$ on $M_0=(0,\infty)\x\pl\Mbar$ has an inverse given 
in terms of eigendecomposition of $\Delta_{\pl\Mbar}$: if $E_{\nu_j}:=\ker(\Delta_{\pl\Mbar}-\nu_j^2+(n-2)^2/4)$ 
where $\nu_j^2\in N_\pl$, that is 
\begin{equation}\label{bfo-conic-model}
\begin{gathered}
 Q_{\bfo}:= \sum_{j=0}^\infty\Pi_{E_j}\Big(I_{\nu_j}(\kappa)K_{\nu_j}(\kappa')H(\kappa'-\kappa)+
I_{\nu_j}(\kappa')K_{\nu_j}(\kappa)H(\kappa-\kappa')\Big)
\\ \times \left|\frac{d\kappa dyd\kappa'dy'}{\kappa\kappa'}\right|^\demi
\end{gathered}\end{equation}
solves $P_{\bfo} Q_{\bfo} = \Id$ on b-half densities on $M_0\x M_0$. 
The index set of $Q_{\bfo}$ at $\zf$ is of the form $0+\mc{Z}$ where $\mc{Z}$ is an index set 
satisfying $\mc{Z}\geq \min(2,n-2)-\eps$ for any $\eps>0$.
Therefore we set
\begin{equation*}
G_{\bfo}^{-2} = (\kappa \kappa') Q_{\bfo} .
\end{equation*}

\begin{remark}\label{many-ends} If $M^\circ$ has $j \geq 2$ ends, i.e. $\partial M$ has $j \geq 2$ components $S_1, \dots, S_j$, then $\bfo$ has $j^2$ components. In this case, all of the off-diagonal components of $G_{\bfo}^{-2}$ are zero, while the diagonal components are completely independent: the $i$th diagonal component only depends on the metric $h(0)$ at $S_i$. 
\end{remark}

\subsection{Leading term at $\rbo$ and $\lbo$.}\label{rbolbo} The leading terms  $G^{n/2-2}_{\rbo}$ and $G^{n/2-2}_{\lbo}$ 
are defined similarly to the previous section since $\nu_0=(n-2)/2$ here as in the asymptotically Euclidean case. However, we remark on what happens when there is more than one end, say $j \geq 2$ ends. Then from matching with $G_{\bfo}^{-2}$ we obtain  boundary conditions for the function $v(z)$ from Section~\ref{34} on the $i$th component of $\rbo$:  namely, it should tend to a nonzero limit at $S_i$ and to zero at all other ends (cf. Remark~\ref{many-ends}). Again Theorem~\ref{rit} implies that there is a unique function (up to scaling) with this property. 

This is sufficient
to deal with dimension $n\geq 4$ since we gain an error that has a positive order of vanishing at $\rbo,\lbo$.
However, to deal with dimension $n=3$, it is not sufficient and we have to construct the terms 
$G_{\rbo}^{\nu_j-1}$ such that $0<\nu_j\leq 1/2$ in the parametrix. Actually it is possible 
to compute the actual second asymptotic term at $\rbo$ in any dimension by this method, when $\nu_1-\nu_0<1$, 
so we will write it for completeness and because it will imply an interesting result for Riesz transform.\\  

For that, we notice that $Q_b$ has asymptotic expansion near $\rbo$ 
\[Q_b=\sum_{\nu_0\leq \nu_j<\nu_0+1}{x'}^{\nu_j}v_{\nu_j}(z,y')+O({x'}^{\nu_0+1+\eps})\]
for some $\eps>0$ and $P_bv_{\nu_j}(z,y')=0$ with $v_{\nu_j}$ growing like $x^{-\nu_j}$.
Checking consistency at $\bfo$ with $I_{\rm ff}(Q_b)$ shows that 
\begin{equation}\label{vnuj}
v_{\nu_j}(x,y,y')= \frac{\Pi_{E_j}(y,y')}{2\nu_j}x^{-\nu_j}+O(x^{-\nu_{j-1}}\log x). 
\end{equation}
Now it suffices to set 
\[G_{\rbo}^{\nu_j-1}:= {x}^{-1}\frac{\kappa'K_{\nu_j}(\kappa')}{\Gamma(\nu_j)2^{\nu_j-1}}v_{\nu_j}(z,y,y')\]
which matches with $G_{\bfo}^{-2}$ and $G_{\zf}^0$.
 
\subsection{Consistency} 
This is almost the same as in Section~\ref{nokernel}; we just have to change 
the normal operator of $Q_b$ on $\zf\cap\bfo$
since the normal operator of $P_b$ has the expression $\oplus_{j\in\nn_0}((sD_s)^2 + \nu_j^2)$ (where $s:=\beta^*(x/x')\in\rr$ is a coordinate on the fibers of ff) after decomposing according to the eigendecomposition of $\Delta_{\pl\Mbar}$. The solution is given by standard Sturm-Liouville theory 
\begin{equation*}
\sum_{j=0}^{\infty} \Pi_{E_j} \frac{e^{-\nu_j |\log s|}}{2\nu_j}|dydy'ds/s|^{\demi}.
\end{equation*}

Then everything is similar to the case $\pl\Mbar=S^{n-1}$ except that $G^{-2}_{\bfo}$ and $G^0_{\zf}$ are polyhomogeneous at $\lbo$ and $\rbo$ with another index set $\mc{E}_{\rb}(Q_b)-1\geq n/2-2$ which does not necessarily contain half integers. 

Therefore, there exists a $G = G(k) \in \Psi^{-2, (-2, 0, 0); \mc{G}}(M, \Omegab)$ satisfying $\mc{G}_{\sca} =0$, $\mc{G}_{\zf} \subset 0\cup \mc{Z}$, $\mc{G}_{\bfo} = -2$, $\mc{G}_{\lbo} = \mc{G}_{\rbo} \subset \cup_{\nu_0\leq\nu_j<\nu_0+1}(\nu_j,0) \cup \nu_0+\mc{Z}$ for some index set $\mc{Z}>-\eps$ for any $\eps>0$, and $G(k)$ with the values at each face prescribed by our models.

\subsection{Exact resolvent}\label{exact-resolvent-conic}
Applying $P + k^2$ to $G(k)$, we define the error $E(k)$ by
$$
(P + k^2) G(k) = \Id + E(k)
$$
where $E(k) \in \Psi^{-\infty, (1, 1, 1); \mc{E}}(M, \Omegab)$, and $\mc{E}$ satisfies $\mc{E}_{\zf} \geq 1-\eps$ for any $\eps>0$, $\mc{E}_{\bfo} \geq 1$, $\mc{E}_{\sca} \geq 1$, $\mc{E}_{\lbo} \geq n/2$, $\mc{E}_{\rbo} \geq n/2-1-\eps$ 
for any $\eps>0$. 
In particular, 
$\mc{E}_{\textrm{f}} > 0$ for each boundary hypersurface $f$ and  Corollary~\ref{cor} shows that the method  detailed for $\pl\Mbar=S^{n-1}$ applies, and the resolvent $R(k)$ has the same properties except the index set which becomes  
\begin{equation}\label{rif2}
\begin{gathered}
\mc{R}_{\zf}  \subset (0,0) \cup (1+\mc{F}), \quad 
\mc{R}_{\lbo},\mc{R}_{\rbo}  \subset \bigcup_{\nu_0\leq\nu_j<\nu_{0}+1}(\nu_j-1) \cup (n/2-1+\mc{F}),\\
\mc{R}_{\bfo}  \subset (-2,0)\cup (-1+\mc{F}), \quad
\mc{R}_{\sca} = 0
\end{gathered}
\end{equation}
for some index set $\mc{F}\geq -\eps$  for any $\eps>0$. 
Moreover  Proposition \ref{prop:composition} shows that the first asymptotic terms $R_{\rbo}^{\nu_j-1}$ 
of the resolvent at $\rbo$ are given by $G_{\rbo}^{\nu_j-1}$ for $\nu_j<\nu_0+1=n/2-1$.

\begin{remark}\label{exact-cones} Suppose that $(M^\circ, g)$ is a metric cone. Let $\Delta_0$ be the Laplacian on $M^\circ$, defined as the Friedrichs extension of the operator with domain given by smooth functions with support compact and disjoint from the cone point. (This operator is essentially self-adjoint in dimensions $n \geq 4$ but not for $n =2,3$.) The resolvent $R_0 = (\Delta_0 + k^2)^{-1}$ scales as $\lambda^{-2}$ under the transformation $x \to \lambda x$, $x' \to \lambda x'$ and $k \to \lambda k$ and is therefore determined by its leading term at $\bfo$, $(R_0)_{\bfo}^{-2}$. It is not hard to see, by expanding in eigenfunctions on the cross-section $\partial M$,  that this model is given precisely by \eqref{bfo-conic-model}. Now choose a cutoff function $\phi$ that is equal to $1$ near infinity and $0$ near the cone point. Then passing to $\phi R_0 \phi$ removes singularities of the kernel of $R_0$ due to the cone point, and we see that this kernel satisfies \emph{all the properties of  the resolvent listed above}. 
We shall use this fact in the following section. 
\end{remark}



\section{Riesz transform}

\subsection{Boundedness on $L^p$ of a class of operators}

We first determine a class of operators in our calculus which are bounded on $L^p(M)$. 
\begin{prop}\label{kint} Let $\alpha, \beta,\delta > 0$. 
Assume that $A = A(k) \in \Psi_{k,\sca}^{-1, (-1, 0, 0), \mc{A}}(M, \Omegabht)$ with index family $\mc{A}$ satisfying
$$
\mc{A}_{\lb} = \mc{A}_{\rb} = \mc{A}_{\bfc} = \emptyset, \quad \mc{A}_{\bfo} > -1, \ \mc{A}_{\sca} \geq 0, \ \mc{A}_{\zf} >-1 , $$ $$ \mc{A}_{\rbo} \geq \frac{n}{2} - 1 - \alpha, \ \mc{A}_{\lbo} \geq \frac{n}{2} -1 - \beta . 
$$
Then $\int_0^{k_0} A(k) \, dk$ is bounded on $L^p(M)$ for
$$
\frac{n}{n-\beta} < p < \frac{n}{\alpha}.
$$
\end{prop}

\begin{proof}
Let us first set $\delta>0$ so that $\mc{A}_{\zf}>\delta$ and $\mc{A}_{\bfo}>\delta$.
We give the details when the kernel is localized to the region $x \leq \epsilon, x' \leq \epsilon$. 
When the kernel of $A$ is localized to the region $x \geq \epsilon/2,  x' \geq \epsilon/2$ then $A(k)$ is a bounded family of properly supported pseudodifferential operators of order $-1$ and the result is classical, while if we localize the kernel to $x \geq \epsilon$, $x' \leq \epsilon/2$ or to 
$x' \geq \epsilon$, $x \leq \epsilon/2$ then the argument is similar (but simpler) to that below for the kernels $A_1$ and $A_2$. 

In the region $x \leq \epsilon$, $x' \leq \epsilon$ we can use coordinates $(x,y)$ in a collar neighbourhood of the boundary of $M$. We define $s = x/x'$ and write $dy$ for the measure density on $\pl\Mbar$.

Let us decompose $A = A_\Delta + A_{\phg}$ where $A_\Delta$ is supported close to the diagonal $\Delta_{k,\sca} \subset \MMksc$ and $A_{\phg}$ is polyhomogeneous on $\MMksc$ with the same index family as $A$. As above we assume both kernels are supported in $x \leq \epsilon, x' < \epsilon$. We further decompose $A_{\phg} = A_1 + A_2 + A_3$, where $A_1$ is supported away from $\sca$ and in $s \geq 1$, $A_2$ is supported away from $\sca$ and in $s \leq 1$ and $A_3$ is localized near $\sca$. We decompose $A_\Delta = A_4 + A_5$ where $A_4$ is localized in $k > x'$ and $A_5$ in $k < x'$. Note that we can assume that $1/2 < s < 2$ on the support of $A_4$ and $A_5$. 

It is natural here to switch back to writing the kernel of $A$ as a multiple of the scattering half-density $|dg dg'|^{1/2}$, since our $L^p(M)$ spaces are with respect to the Riemannian measure $|dg|$. Indeed we ignore the half-density factors from here on (which are  not natural for $p \neq 2$). Then the order of vanishing of the kernel is $n-1 + \epsilon$ at $\bfo$ for some $\epsilon > 0$, $n -1-\alpha$ at $\rbo$ and $n-1-\beta$ at $\lbo$, and $0$ at $\sca$ and $\zf$, due to the relation $|dg_b dg_b'|^{1/2} = (xx')^{n/2} |dg dg'|^{1/2} $.  

Consider the kernel $A_1$. Since it is supported away from $\sca$ we may use boundary defining functions  $x'/(x'+k)$ for $\rb$, $x+x'+k$ for $\bfo$ and $(x'+k)/(x+x'+k)$ for $\rbo$. So we can bound the kernel of $\int A_1 dk$ by 
\begin{equation*}\begin{gathered}
 C_N\int_{0}^{k_0}\Big(\frac{k}{x'+k}\Big)^{-1+\delta}(\frac{x'}{x'+k})^N (x+x'+k)^{n-1+\delta}
\Big(\frac{x'+k}{x+x'+k}\Big)^{n-1-\alpha} \, dk \\
\leq C_N  {x'}^{n+\delta} s^{\alpha+\delta}\int_{0}^\infty \bar{k}^{-1+\delta}(1+\overline{k})^{-N+n} \, d\overline{k}\leq 
C_N  {x'}^{n+\delta} s^{\alpha+\delta}
\end{gathered}\end{equation*}
where $N>0$ is arbitrary and we used change of variable $\overline{k}x'= k$.
We look at $L^p$ boundedness of $\int_{0}^kA_1(k)dk$. Let $f(x,y)$ be a function on $[0, \epsilon]_x \times \pl\Mbar$. Then we need to show 
\begin{equation*}\label{estimlp}
\int \!\!\! \int_{0}^1 \Bigg|\int \!\!\! \int_{0}^x f(x',y')s^{\alpha+\delta}{x'}^{-1+\delta}dx'dy' \Bigg|^px^{-n-1}dxdy\leq C||f||^p_{L^p(M)}.\end{equation*}
We bound the left hand side using H\"older's equality (with $q = p/(p-1)$) to obtain
\begin{equation*}\begin{gathered}
\int \!\!\! \int_{0}^1 \Bigg|\int  \!\!\! \int_{0}^x f(x',y'){x'}^{-1-\alpha}dx'dy' \Bigg|^px^{-n-1+p(\alpha+\delta)}dxdy  \\
\leq 
\int \!\!\! \int_{0}^1 \Big(\int \!\!\!
\int_{0}^x |f(x',y')|^p \frac{dx' dy'}{  {x'}^{n+1}} \Big) \Big(\int \!\!\! \int_{0}^x {x'}^{q(-1-\alpha+(n+1)/p)} dx' dy' \Big)^{\frac{p}{q}}
 x^{p(\alpha+\delta)}\frac{dxdy}{ x^{n+1}} \\
 \leq C \| f \|_p^p \int \!\!\! \int_{0}^1 x^{p(-1-\alpha)+n+1 + (p-1)} x^{-n-1+p(\alpha+\delta)}dxdy \\
 = C \| f \|_p^p \int \!\!\! \int_{0}^1 x^{p\delta - 1} dxdy  \leq C' 
 \| f \|_p^p.
 \end{gathered}\end{equation*}
where we used $p<n/\alpha$ to integrate the $x'$ integrand of the second line.

The same argument applied to the formal adjoint of $A_2$ shows that $A_2$ is bounded on $L^p$ for $p > n/(n-\beta)$. 


For the kernel $A_3$, we may localize near a fiber of $\sca$ with basis point $y'$, then if $y_1,\dots, y_n$ are local coordinates near $y'$, we can use the boundary defining function $(1 + k^2 |z-z'|^2)^{-1/2}$ for $\bfc$, where $z_0=1/x-1/x'$
and $z_i=y_i/x-y_i/x'$, we also use $x+x'+k$ for boundary defining function of $\bfo$. However, since $A_3$ is by assumption supported near $\sca$, we may assume that $k>(x+x')$ on the support on $A_3$ and hence we may take $k$ as a boundary defining function for $\bfo$. This allows us to bound $\int A_3(k) \, dk$ by
$$
\int_0^{k_0} \big( \frac{1}{1 + k^2 |z-z'|^2} \big)^N k^{n-1+\delta} \, dk
$$
for any $N$. By changing variable to $k|z-z'|$ in the integral we see that this is bounded by $\min(1, |z-z|^{-n-\epsilon})$, which shows that it is bounded on $L^p$ for $1 \leq p \leq \infty$.

For the kernel $A_4$, we use same notations for $z$ and use coordinates $W:=k(z-z'),y,{\kappa'}^{-1}=x'/k,k$ as coordinates. Then $W=0$ defines the diagonal. 
Since $A_4$ is conormal to the diagonal of conormal order $-1$, and vanishes to order $n-1+\delta$ at $\bfo$,  the kernel is bounded by
$$
k^{n-1+\delta} h(|W|) = k^{n-1+\delta} h(k|z-z'|)
$$
where $h(t)$ is bounded by $C t^{-n+1}$ for $t$ small and is rapidly decreasing in $t$ for $t$ large. A similar argument to that for $A_3$ shows that the kernel of $A_4$ is bounded by $\min(|z-z'|^{-n+1}, |z-z'|^{-n-\delta})$, which shows that it is bounded on $L^p$ for $1 \leq p \leq \infty$.

Finally, for $A_5$, localized to $k < x'$,  we 
use coordinates $(s=x/x',x,\kappa'=k/x',y,y')$. Here $y-y'$ and $s-1$, or equivalently $x(z-z')$, define the diagonal  and $x$ is a boundary defining function for $\bfo$ on the support of $A_5$. One can suppose that $A_5$ is supported in $\{x|z-z'|<1\}$. 
The kernel $A_5(k)$ is therefore bounded by
$(k/x)^{-1+\delta}x^{n-1+\delta} h(x|z-z'|)$ where $h$ is as above. Integrating in $k$ from $0$ to $x$ (since $A_5$ is supported where $k \leq x$) gives a bound
$$
\Big| \int_0^x A_5(k) \, dk \Big| \leq x^{n+\delta} h(x|z-z'|).
$$
This is bounded by $\min(|z-z'|^{-n+1}, |z-z'|^{-n-\delta})$, which shows that it is bounded on $L^p$ for $1 \leq p \leq \infty$. This completes the proof of the proposition. 
\end{proof}

\subsection{Boundedness of the Riesz transform} 
Recall from the Introduction that we define the Riesz transform $T$ by 
$T = d \circ (P_>)^{-1/2}$. This can be written, at least formally, as
\begin{equation}
T = \frac{2}{\pi} d \int_0^\infty R(k) \Pi_> \, dk
\label{T}\end{equation}
where $\Pi_>$ is the spectral projection onto $(0, \infty)$ for the operator $P$. 

Using the elementary bound $\| R(k) \Pi_> \|_{L^2 \to L^2} \leq k^{-2}$, we see that the integral converges on $[k_0, \infty)$ for any $k_0 > 0$, at least as an operator $L^2(M) \to H^{-1}(M)$. To make sense of the integral near zero we use the results in this paper on the resolvent kernel $R(k)$ as $k \to 0$. 

To compute the integral \eqref{T} we divide it into pieces and also compare it to the classical Riesz transform. We recall that $P$ has at  most finitely many negative eigenvalues $0 > -k_1^2 \geq -k_2^2 \geq \dots \geq -k_M^2$. 
We write $\Pi_j$ for the projection onto the $-k_j^2$ eigenspace. 
Choose a smooth monotone function $f(t)$ that is equal to $0$ for $t \leq -k_1^2/2$ and $1$ for $t \geq -k_1^2/4$. Thus $f(P) = \Pi_>$. We also choose a smooth function $\psi(k)$ so that $\psi(k) = 1$ for $k \geq k_1^2 /2$, $\psi(k) = 0$ for $k \leq k_1^2/3$.

Consider the integral 
\begin{equation}
\frac{2}{\pi} d  \int_0^\infty \psi(k) R(k) \Pi_> \, dk = 
\frac{2}{\pi} d  \int_0^\infty \psi(k) R(k) f(P) \, dk,
\label{klarge}\end{equation}
where we inserted the cutoff function $\psi$ into integral \eqref{T} for $T$. 
Due to the conditions on $f$, the function
$$
g(\lambda) =  f(\lambda) \int_0^\infty \psi(k) \frac1{k^2 + \lambda^2} \, dk
$$
is a smooth function of $\lambda$ which is a classical symbol of order $-1$ as $\lambda \to \infty$. By the symbolic functional calculus of \cite{HV}, the operator $g(P)$ is a scattering pseudodifferential operator of order $-1$. Hence $d \circ g(P)$,
which is equal to \eqref{klarge}, is a scattering pseudodifferential operator of order $0$. 
Using Stein \cite{Stein}, chapter 6, section 5, this is bounded\footnote{The result in Stein is stated for Euclidean space; however, by localizing in $\partial \Mbar$ one reduces to this case.} 
 on $L^p(M)$ for $1 < p < \infty$. This treats the integral for  \eqref{T} for large $k$.
 
 Now consider the operator $ \phi R_0(k) \phi$, where $R_0(k)$ is the resolvent for the free Laplacian $\Delta_0$ on the exact cone $(0,\infty)\x\pl\Mbar$ and $\phi$ is a cutoff function, equal to zero on $\{x>\eps\}$ and $1$ on $\{x<\eps/2\}$ for some small $\eps>0$. If $\eps$ is small enough, this kernel may be regarded as living on $\MMksc$ and it satisfies the properties listed in Section~\ref{exact-resolvent-conic} (see Remark~\ref{exact-cones}).  We can apply the argument above to $R_0(k)$ and we see that $d \circ g(\Delta_0)$ is 
 bounded on $L^p(M)$ for $1 < p < \infty$.
  
Now consider the remaining integral, where $1 - \psi(k)$ is inserted into the integral \eqref{T} for $T$. Choosing $k_0 = k_1^2/3$ we can write this as
$$
\frac{2}{\pi} d \int_0^{k_0} (1 - \psi(k)) R(k) \Big( \Id  - \sum_{j=1}^M \Pi_j \Big) \, dk.
$$
It is well known via Agmon estimates \cite{Agmon}  that the eigenfunctions corresponding to the $\Pi_j$, $j \geq 1$, together with all their derivatives,  are exponentially decreasing  at infinity. It follows that the $\Pi_j$ in the integral above contribute an operator that is bounded on $L^p$ for $1 \leq p \leq \infty$. The remainder is
$$
\frac{2}{\pi} d \int_0^{k_0} (1 - \psi(k)) R(k)  \, dk.
$$
Rather than study this directly, we subtract off the free resolvent kernel, and consider
\begin{equation}
\frac{2}{\pi}  \int_0^{\infty} \Big( d(R(k)  - \phi R_0(k) \phi) 
+[d,\phi]R_0(k)\phi \Big) dk.
\label{r-diff}\end{equation}
By Theorems ~\ref{v=0} and ~\ref{th2}, $R(k)$ is in the calculus with pseudodifferential order $-2$ and with index sets $\geq 0$ at $\zf$, $\geq n/2 -2$ at $\rbo$ and $\lbo$, $\geq -2$ at $\bfo$, $\geq 0$ at $\sca$ and trivial at all other boundary hypersurfaces. Similarly, our analysis shows that $[d,\phi]R_0(k)\phi$ and $\phi R_0(k) \phi $ are in the calculus with pseudodifferential order $-2$ and with similar index sets, except that $[d,\phi]$ is compactly supported so
$[d,\phi]R_0(k)\phi$ vanishes at all order at all boundary hypersurfaces but $\rb,\rbo,\zf$. The point of subtracting the two resolvents is that they have \emph{the same leading term at $\bfo$ and $\sca$}, since these models are determined purely by the metric $h(0)$. Hence the difference has index set $\geq -1$ at $\bfo$ and $\geq 1$ at $\sca$\footnote{Actually, there could be a log term at order $-1$ at $\bfo$. However, this would make no difference to the argument below, so we ignore this detail to simplify the exposition.}. 

Now consider what happens when $d$ is applied on the left. 
The operator $d$ has components $\partial_{z_i}$ each of which lift from the left factor to be of the form $\rho_{\bfo}\rho_{\lbo}\rho_{\lb}\rho_{\bfa}$ times a vector field tangent to the boundary of $\MMksc$, and transverse to the diagonal. 
Therefore,  it increases the order of vanishing to $n/2 - 1$ at $\lbo$ and to $0$ at $\bfo$. At $\rbo$, by contrast, it acts globally in the $z$ variable. So we get an improvement at $\rbo$ if and only if the leading coefficient is annihilated by $d$, or equivalently is constant in $z$. (It also increases the pseudodifferential order to $-1$). 
The leading term of $d(R(k)-\phi R_0(k)\phi-[d,\phi]R_0(k)\phi)$ at $\rbo$ is at order $n/2-2$, this is 
$ dR_{\rbo}^{\ndemi-2}-\phi(z)dR_{0,\rbo}^{\ndemi-2}$
where $R_{\rbo,0}^{\ndemi-2}$ is the leading term of $R_0(k)$ at $\rbo$.
It is important to notice that this term vanishes if and only if $R_{\rbo}^{\ndemi-2}$ and $R_{0,\rbo}^{\ndemi-2}$ are constant in the variable $z$, which is exactly the case when $V\equiv 0$ and $M$ has one end.
Indeed, we have seen in Sections~\ref{34} and \ref{rbolbo} that $R_{\rbo}^{\ndemi - 2}$, as a function of $z$, is a multiple of a bounded function annihilated by $P$. This is uniquely determined up to constants when $M$ has one end, and it is constant if and only if $V \equiv 0$ (see Remark~\ref{vzconstant}). When $M$ has more than one end, on the $i$th component of $\rbo$ it is a multiple of the bounded solution of $P v = 0$ that tends to $1$ at the $i$th end and $0$ at all other ends (see Remark~\ref{many-ends}). In particular, in this case it is not constant.
When it is constant, it means from Subsection \ref{rbolbo} that the next leading term at $\rbo$ is at order $\nu_1-1$, 
so we can now apply Proposition~\ref{kint} to deduce that \eqref{r-diff} is bounded for $1<p<n$ when $\nu_1\geq 1$
and for $1<p<n/(n/2-\nu_1) := p_{\max}$ otherwise. 
Combined with the earlier results about the integral for large $k$, we have proved that $T - \phi T_0 \phi$ is bounded on $L^p$ for $1<p<p_{\max}$, where $T_0$ is the Riesz transform on the exact cone $\rr^+\x\pl\Mbar$. When $M$ is asymptotically Euclidean, the cone is Euclidean $\RR^n$ and the boundedness of $T_0$ for $1 < p < \infty$ is classical. For the general conic case, it is shown in  \cite{Li} that $T_0$ is bounded on $L^p$ for $1 < p < p_{\max}$  so we conclude that $T$ itself is bounded on $L^p$ for the stated range.\\

Finally we claim that the range of $p$ is sharp. In this case, we can write the resolvent kernel near the right boundary (now as a multiple of the scattering half-density rather than the b-half density, which gives an extra factor of $(xx')^{n/2}$) as 
\begin{equation*}\begin{gathered}
R(k) =  k^{\ndemi-2}\frac{  \kappa' K_{n/2 - 1}(\kappa')}{\Gamma(n/2 - 1) 2^{n/2 -2}}{x'}^{\ndemi}x^{\ndemi-1} v_{\nu_0}(z)+
k^{\nu_1-1}\frac{\kappa' K_{\nu_1}(\kappa')}{\Gamma(\nu_1) 2^{\nu_1-1}}{x'}^{\ndemi}x^{\ndemi-1} v_{\nu_1}(z)
\\ = {x'}^{n-2}\frac{{\kappa'}^{\ndemi-1} K_{\ndemi-1}(\kappa')}{\Gamma(\ndemi-1)2^{\ndemi-2}}x^{n/2-1}v_{\nu_0}(z) +
{x'}^{\ndemi+\nu_1-1}\frac{{\kappa'}^{\nu_1} K_{\nu_1}(\kappa')}{\Gamma(\nu_1)2^{\nu_1-1}}x^{n/2-1}v_{\nu_1}(z) 
\end{gathered}\end{equation*}
modulo $O(\rho_{\rbo}^{n-2+\nu_2})+O(\rho_{\rbo}^{n-1}\log \rho_{\rbo})$.  So we have 
\begin{multline*}
d R(k)   = {x'}^{n-2}\frac{{\kappa'}^{\ndemi-1} K_{\ndemi-1}(\kappa')}{\Gamma(\ndemi-1)2^{\ndemi-2}}d(x^{n/2-1}v_{
\nu_0}(z))\\
+{x'}^{\ndemi+\nu_1-1}\frac{{\kappa'}^{\nu_1} K_{\nu_1}(\kappa')}{\Gamma(\nu_1)2^{\nu_1-1}}d(x^{n/2-1}v_{\nu_1}(z))
\end{multline*}
(since $d$ acts in the left variable $z$ here) modulo $O(\rho_{\rbo}^{\ndemi+\nu_2})+O(\rho_{\rbo}^{n-1}\log \rho_{\rbo})$. Integrating in $k$, we change variable to $\kappa' = k/x'$ (which gives us an extra power of $x'$) and obtain 
\begin{multline*}
 {x'}^{n-1} d(x^{n/2-1}v_{
\nu_0}(z))\int_0^\infty  \frac{{\kappa'}^{\ndemi-1} K_{\ndemi-1}(\kappa')}{\Gamma(\ndemi-1)2^{\ndemi-2}}\, d\kappa'\\  
+{x'}^{\ndemi+\nu_1}d(x^{n/2-1}v_{
\nu_1}(z))\int_0^\infty  \frac{{\kappa'}^{\nu_1} K_{\nu_1}(\kappa')}{\Gamma(\nu_1)2^{\nu_1-1}}\, d\kappa'
\end{multline*}
modulo $O({x'}^{\ndemi+\nu_2})+O({x'}^{n}\log x')$ for the asymptotics of $T$ at $x' = 0$. 
Since $x^{n/2 - 1} v_{\nu_0}(z)$ is constant if and only if $V\equiv 0$ and $M$ has one end, 
$d (x^{n/2 - 1} \varphi_j(z))$ is identically zero if only and if we are in this case. 
In addition, the $\kappa'$ integrals do not vanish, indeed
the function (for $\nu>0$)
\[F(\kappa'):={\kappa'}^{\nu}K_{\nu}(\kappa')\]
is smooth on $(0,\infty)$, with $F(0)=2^{\nu-1}\Gamma(\nu)$ 
and $F(\kappa)=F(0)+O(\kappa^2)+O(\kappa^{2\nu})$ as $\kappa\to 0$, and we have the identity for $\nu>0$
\[\int_{0}^\infty F(\kappa')d\kappa'=-2^{\nu-1}\pi^\demi\Gamma(\nu+\demi)\]
since $z^{\nu}K_{\nu}(z)=c_\nu\mc{F}_{t\to z}((1+t^2)^{-\nu-\demi})$ with $c_\nu:=-\Gamma(\nu+\demi)2^{\nu-1}\pi^{-\demi}$.

Hence $T = a(z, y') {x'}^{n-1} + O({x'}^{\nu_1+\ndemi})+O({x'}^{n}\log x')$ at $x' = 0$ where $a$ does not vanish
if and only if $V\not\equiv 0$ or $M$ has more than one end.
It follows immediately that the upper threshold for $p$ is sharp in that case. A similar analysis at the left boundary shows that $T = b(y, z') x^{n-2+m'} + O(x^{n-1+m'})$ as $x \to 0$ where $b$ does not vanish, showing that the lower threshold is also sharp. 

Now if $V\equiv 0$ and $M$ has one end, $T=c(z,y'){x'}^{\nu_1+\ndemi}+O({x'}^{\nu_2+\ndemi})+O({x'}^{n-1}\log x')$
with $c(z,y')=d_{z}(x^{\ndemi-1}v_{\nu_2}(z,y'))$ and $v_{\nu_2}$ solving $P(x^{\ndemi-1}v_{\nu_2}(z,y'))=0$ and $v_{\nu_2}(x,y,y')\sim \gamma(y',y)x^{-\nu_2}$ as $x\to 0$ for some $\gamma(y,y')\not=0$. Thus $c\not\equiv 0$ and this proves the upper range is sharp in that case, at least when $\nu_1<\nu_0+1=\ndemi$. 

This completes the proof of the theorem.\\


\section{Appendix: proof of Proposition~\ref{prop:composition}}\label{comp-proof}

This result will follow from the existence of a `triple space' $M^3_{k, \sca}$ having the property that there are b-fibrations down to $M^2_{k, \sca}$ which are lifts of the three projections $M^3 \times [0, k_0] \to M^2 \times [0, k_0]$. We now proceed to define such a space.

We first define $M^3_{k, b}$ to be the `total blowup' of the space $M^3 \times [0, k_0]$ in the sense of \cite{asatet}, that is, with the codimension 4 corner blown up, followed by the 4 codimension 3 corners (which may be done in any order as they are separated after the first blowup), followed by the 6 codimension 2 corners (which likewise may be done in any order). Then the three projections from $M^3 \times [0, k_0] \to M^2 \times [0, k_0]$ extend from the interior to b-fibrations from  $M^3_{k,b} \to M^2_{k, b}$ (this can be deduced from the lemmas in section 2.4 of \cite{asatet}, for example). Let us denote these b-fibrations $\pi_{b, *}$ for $* = L, C, R$ (left, centre and right) according as it eliminates the right, centre or left variable respectively. 

The space $\MMksc$ is the space $M^2_{k,b}$ with the boundary $\partial \Delta_{k,b}$ at $\bfc$ blown up. The inverse image of this submanifold in $M^3_{k,b}$ under any one b-fibration $\pi_{b, *}$ is the union of 2 p-submanifolds, $J_*$ and $G_*$. Their intersection properties are such that once the intersection $K$ of the $G_*$ is blown up, the lifts of $G_*$ are all separated, while the $J_*$ only intersect the corresponding $G_*$. Therefore, we define
\begin{equation}
M^3_{k, \sca} = \big[ M^3_{k,b}; K; G_L, G_C, G_R; J_L, J_C, J_R \big].
\label{triple-blowup}\end{equation}

\begin{lem}\label{lem:bfib} The three projections from the interior of  $M^3 \times [0, k_0]$ to the interior of $M^2 \times [0, k_0]$ extend by continuity to b-fibrations from $M^3_{k, \sca}$ to $\MMksc$. 
\end{lem}

The proof uses the following properties of blowups and b-fibrations (see \cite{scatmet} or \cite{asatet}). 

\begin{lem}\label{sss} Suppose that $X$ and $Y$ are manifolds with corners and $f : X \to Y$ is a b-fibration. Let $S$ and $T$ be p-submanifolds of $X$. 
If either 

(i) $S \subset T$, or 

(ii) $S$ is transverse to $T$, 

then $[X; S; T] = [X; T; S]$.

(iii) Suppose that $S$ is b-transversal to $f$ and such that $f(S)$ is not contained in any codimension 2 face of $Y$. Then $f$ restricted to $X \setminus S$ extends to a b-fibration from $[X;S] \to Y$. 
\end{lem}

\begin{proof}[Proof of Lemma~\ref{lem:bfib}]
The proof runs along the same lines as the proof of the existence of b-fibrations from the triple scattering space $M^3_{\sca}$ to $M^2_{\sca}$ in Section 23 of \cite{scatmet}, and we follow this proof to the extent possible. We refer to \cite{cocdmc} for definitions of terms such as b-fibration and b-transversal. 

Since the spaces $M^3_{k, b}$ and $M^3_{k, \sca}$ are symmetric under permutation of the $M$ factors, it is enough to show that one of the projections from the interior of $M^3 \times [0, k_0]$ to the  interior of $M^2 \times [0, k_0]$ extends to be a b-fibration. We shall take the projection $\pi_R$  that projects off the right factor of $M$.

We first show consider the space $M^3_{k,b}$. We temporarily use 4-digit binary codes for the faces of $M^3 \times [0, k_0]$ where $0000$ represents $(\partial M)^3 \times \{ 0 \}$, $1010$ stands for $M \times \partial M \times M \times \{ 0 \}$, $0011$ stands for $(\partial M)^2 \times M \times [0, k_0]$, etc. Then 
$$
M^3_{k,b} = \big[ M^3 \times [0, k_0] ; 0000; 1000, 0100, 0010, 0001; 1100, 1010, 1001, 0110, 0101, 0011 \big].
$$
Here the four codimension 3 faces may be permuted among themselves since they become disjoint after the $0000$ blowup, and similarly the six codimension 2 faces may be permuted among themselves. So we can write
$$
M^3_{k,b} = \big[ M^3 \times [0, k_0] ; 0000; 0010, 0001, 0100, 1000;  1010,  0110,  0011, 1100, 1001, 0101 \big].
$$
Now the $0000$ and $0010$ blowups may be commuted using (i) of Lemma~\ref{sss}. Then the $1010$ blowup may be moved to the left past $0001$ (disjoint, once $0000$ has been blown up), $0100$ (disjoint), $1000$ (contains) and $0000$ (contains). Similarly $0110$ and $0011$ may be moved to the left, obtaining
$$
M^3_{k,b} = \big[ M^3 \times [0, k_0] ; 0010; 1010,  0110,  0011; 0000;  0001,  0100, 1000;  1100, 1001, 0101 \big].
$$
The first four blowups here give us $M^2_{k,b} \times M$. So we have shown that 
$$
M^3_{k,b} = \big[ M^2_{k,b} \times M;  0000;  0001,  0100, 1000;  1100, 1001, 0101 \big].
$$
$$
 = \big[ M^2_{k,b} \times M;  \bfo \times \partial M;  \bfc \times \partial M,  \lbo \times \partial M,  \rbo \times \partial M;  \zf \times \partial M, \rb \times \partial M, \lb \times \partial M \big].
$$

Now we include the blowups that turn $M^3_{k,b}$ into $M^3_{\sca, b}$. Note that these may be reordered $G_R, J_R, K, G_C, J_C, G_L, J_L$ since $K \subset G_R$ and $K$ and $J_R$ are disjoint. So we have
\begin{multline*}
M^3_{k, \sca} = \big[  M^2_{k,b} \times M; \bfo \times \partial M;  \bfc \times \partial M,  \lbo \times \partial M,  \rbo \times \partial M;  \\ \zf \times \partial M, \rb \times \partial M, \lb \times \partial M; 
G_R; J_R; K; G_C, J_C; G_L, J_L \big].
\end{multline*}
Recall that $G_R$ is the submanifold $\partial_{\bfc} \Delta_{k,b} \times \partial M \subset M^2_{k,b} \times M$, and 
$J_R$ is the submanifold $\partial_{\bfc} \Delta_{k,b} \times M \subset M^2_{k,b} \times M$, lifted via the blowups listed to their left. Since $\partial_{\bfc} \Delta_{k,b} \subset \MMksc$ does not intersect $\zf$, $\lbo$, $\rbo$, $\lb$, $\rb$,  they can be commuted to the left past five blowups. So we get
\begin{multline*}
M^3_{k, \sca} = \big[  M^2_{k,b} \times M; \bfo \times \partial M;  \bfc \times \partial M; G_R; J_R; \lbo \times \partial M,  \rbo \times \partial M; \\  \zf \times \partial M, \rb \times \partial M, \lb \times \partial M; 
 K; G_C, J_C; G_L, J_L \big].
\end{multline*}

We can commute $G_R$ to the left of $\bfc \times \partial M$ since it is contained in $\bfc \times \partial M$; then we can commute $J_R$ to the left of $\bfc \times \partial M$ since these two submanifolds become disjoint after $G_R$ is blown up; and finally we can commute $J_R$ past $G_R$ since it contains $G_R$. Finally, $J_R$ and $\bfo \times \partial M$ are transverse so these two blowups can be commuted, using (ii) of Lemma~\ref{sss}\footnote{This step is the only essential addition to the arguments of \cite{scatmet}}. 
We end up with
\begin{multline*}
M^3_{k, \sca} = \big[  M^2_{k,b} \times M;  \bfo \times \partial M; J_R; G_R;   \bfc \times \partial M;  \lbo \times \partial M,  \rbo \times \partial M; \\  \zf \times \partial M, \rb \times \partial M, \lb \times \partial M; 
 K; G_C, J_C; G_L, J_L \big] \\
 = \big[  M^2_{k,\sca} \times M;  \bfo \times \partial M; G_R;   \bfc \times \partial M;  \lbo \times \partial M,  \rbo \times \partial M; \\  \zf \times \partial M, \rb \times \partial M, \lb \times \partial M; 
 K; G_C, J_C; G_L, J_L \big].
\end{multline*}
Now using (iii) of Lemma~\ref{sss} repeatedly we see that there is a b-fibration from $M^3_{k, \sca}$ to $M^2_{k, \sca}$. This completes the proof. 
\end{proof}

We resume the proof of Proposition~\ref{prop:composition}. We only prove the result in the case that  $(a_{\sca}, a_{\bfo}, a_{\zf}) = (0,0,0)$ since the general case follows easily from this. Composition may be defined in terms of pullbacks and pushforwards. We first need to take care of the density bundles. We have an isomorphism
\begin{equation}
\pi_L^* \Omegab(\MMksc) \otimes \pi_C^* \Omegab(\MMksc) \otimes \pi_R^* \Omegab(\MMksc) =  \tilde \Omega_b ( \MMMksc) \otimes \big| \frac{dk}{k} \big|^{1/2}
\label{iso}\end{equation}
where $\tilde \Omega_b(\MMMksc)$ is the lift of $\Omega_b(M^3_{k,b}$ to $\MMMksc$. This  allows us to express the composition $A \circ B$ as
$$
(A \circ B) \nu = (\pi_{C})_* \Big( \pi_L^* A \otimes \pi_R^* B \otimes \pi_C^* \nu \otimes \big| \frac{dk}{k} \big|^{-1/2} \Big)
$$
where $\nu$ is a b-half-density on $\MMksc$. 

Given $A$ and $B$ as in the proposition, decompose them into $A = A_{\Delta} + A_{\phg}$, $B = B_{\Delta} + B_{\phg}$ where $A_{\Delta}$ and $B_{\Delta}$ are supported close to $\Delta_{k, sc}$ and classical conormal to it, uniformly to the boundary, while $A_{\phg}$ and $B_{\phg}$ are smooth in the interior and polyhomogeneous at the boundary with the given index sets. The composition result for $A_{\phg}$ and $B_{\phg}$ follows directly from the b-fibration property of Lemma~\ref{sss} and the Pushforward theorem (Theorem 5 of \cite{cocdmc}). For example, the second line of \eqref{comp-if} follows since there are two boundary hypersurfaces of $M^3_{k, \sca}$ that project to $\zf$ under $\pi_C$, one of which projects to $\zf$ under both $\pi_L$ and $\pi_R$, and the other of which projects to $\rbo$ under $\pi_L$ and $\lbo$ under $\pi_R$.

Next, the result for $A_{\Delta}$ and $B_{\Delta}$ follows from the geometry of the lifts of $\Delta_{k, sc} \subset \MMksc$ to $M^3_{k, \sca}$ by each of the 3 b-fibrations. Namely, each of the lifts is an interior p-submanifold, each pair of which intersect transversally at the lifted triple diagonal which projects diffeomorphically to $\Delta_{k, \sca}$ under $\pi_C$.

Finally, consider the composition of $A_{\Delta}$ with $B_{\phg}$ (the case $A_{\phg}$ with $B_{\Delta}$ is equivalent). We observe that the composition property holds if $A$ is a differential operator of any order. We can write $A_{\Delta} = A' \circ D  + A''$, where $D$ is an differential operator of order  $2N + m$ which is elliptic as a conormal distribution to $\Delta_{k, \sca}$, $A'$ is supported close to $\Delta_{k, \sca}$ and is conormal of very negative order $-N$, and $A''$ is smooth and supported close to $\Delta_{k, \sca}$. The composition of $A''$ with $B_{\phg}$ can be treated using the Pushforward theorem. The differential operator $D$ maps $B_{\phg}$ to another polyhomogeneous kernel with the same index family. On the other hand, $A'$ has a kernel which is $C^{N-n-1}$. Therefore $A' \circ D \circ B_{\phg}$ is polyhomogeneous to some finite but large order $\sim N$. Since $N$ is arbitrary, this shows polyhomogeneity of $A_{\Delta}  \circ B_{\phg}$.

\end{document}